\pdfoutput=1
\RequirePackage{ifpdf}
\ifpdf % We are running pdfTeX in pdf mode
\documentclass[pdftex]{sigma}
\else
\documentclass{sigma}
\fi

\usepackage{enumerate,cancel}
\usepackage[all]{xy}

\numberwithin{equation}{section}

\newcommand{\cc}[1]{\overline{{#1}}}
\begin{document}

\allowdisplaybreaks

\renewcommand{\PaperNumber}{032}

\FirstPageHeading

\ShortArticleName{On Orbifold Criteria for Symplectic Toric Quotients}

\ArticleName{On Orbifold Criteria for Symplectic Toric Quotients}

\Author{Carla FARSI~$^\dag$, Hans-Christian HERBIG~$^\ddag$ and Christopher SEATON~$^\S$}

\AuthorNameForHeading{C.~Farsi, H.-C.~Herbig and C.~Seaton}

\Address{$^\dag$~Department of Mathematics, University of Colorado at Boulder,
\\
\hphantom{$^\dag$}~Campus Box 395, Boulder, CO 80309-0395, USA}
\EmailD{\href{mailto:farsi@euclid.colorado.edu}{farsi@euclid.colorado.edu}}
\URLaddressD{\url{http://www.colorado.edu/math/people/professors/farsi.html}}

\Address{$^\ddag$~Centre for Quantum Geometry of Moduli Spaces,
\\
\hphantom{$^\ddag$}~Ny Munkegade 118 Building 1530, 8000 Aarhus C, Denmark}
\EmailD{\href{mailto:herbig@imf.au.dk}{herbig@imf.au.dk}}

\Address{$^\S$~Department of Mathematics and Computer Science, Rhodes College,
\\
\hphantom{$^\S$}~2000 N.~Parkway, Memphis, TN 38112, USA} 
\EmailD{\href{mailto:seatonc@rhodes.edu}{seatonc@rhodes.edu}}
\URLaddressD{\url{http://faculty.rhodes.edu/seaton/}}

\ArticleDates{Received August 07, 2012, in f\/inal form April 02, 2013; Published online April 12, 2013}

\Abstract{We introduce the notion of regular symplectomorphism and graded regular symplectomorphism between singular
phase spaces.
Our main concern is to exhibit examples of unitary torus representations whose symplectic quotients cannot be graded
regularly symplectomorphic to the quotient of a~symplectic representation of a~f\/inite group, while the corresponding
GIT quotients are smooth.
Additionally, we relate the question of simplicialness of a~torus representation to Gaussian elimination.}


\Classification{53D20; 58A40; 13A50; 14L24; 57R18}

\vspace{-2mm}

\section{Introduction}
%\label{intro}

Let $G$ be a~compact Lie group acting on a~symplectic manifold $(M,\omega)$ by symplectomorphisms.
One says that the action is \emph{Hamiltonian} with \emph{moment map} $J\colon  M\to\mathfrak g^*$, $\mathfrak
g^*$ being the dual space of the Lie algebra $\mathfrak g$ of $G$, if
\begin{enumerate}\itemsep=0pt
\item $J$ is a~smooth $G$-equivariant map, \item For each $\xi\in\mathfrak g$ the vector f\/ield
$\{J_\xi,\:\}$ coincides with the fundamental vector f\/ield of $\xi$ acting on $M$, where $J_\xi:=\langle
J,\xi\rangle\in\mathcal C^\infty(M)$ and $\{\:,\:\}$ denotes the Poisson bracket associated to the symplectic form
$\omega$.
\end{enumerate}
The \emph{symplectic quotient} $M_0=Z/G$ is def\/ined to be the space of $G$-orbits in the zero f\/ibre $Z:=J^{-1}(0)$
of the moment map.

It is well-known~\cite{MarsdenWeinstein, Meyer} that if $0\in\mathfrak g^*$ is a~regular value of $J$, then the quotient
$M_0=Z/G$ of the closed submanifold $Z$ by the action of $G$ is in a~canonical way a~symplectic orbifold.
This is the case, for instance, when the $G$-action is locally free.
If $0\in\mathfrak g^*$ is not a~regular value, a~theorem of E.~Lerman and R.~Sjamaar~\cite{LermanSjamaar} tells us that
$M_0=Z/G$ is a~stratif\/ied symplectic space; for more details see Subsections~\ref{reptheory} and~\ref{PoissonDF}.
Note that $0\in\mathfrak g^*$ is a~singular value if, for example, $(M,\omega)$ is a~symplectic vector space, the
$G$-action is linear, and the moment map is chosen to be homogeneous quadratic.
We refer to this situation as \emph{the linear case}.

It has been observed that in the linear case, the symplectic quotient can occasionally be identif\/ied,
symplectically~\cite{CushmanSjamaar,BosGotay,LermanSjamaarMontgomery} or merely
topologically~\cite{HerbigIyengarPflaum}, with a~quotient by a~symplectic representation of a~f\/inite group.
This is the case, for instance, with the physically interesting example of angular momentum~\cite{BosGotay}.
For more examples, see Subsection~\ref{magic}.

Our paper is an attempt towards a~more systematic understanding of when and how this happens.
If one is searching for orbifold criteria, a~natural idea is to use intuition from complex algebraic toric geometry
(see e.g.~\cite{CoxLittleSchenck, FultonBook}).
Namely, if one considers a~representation of a~complex torus $\mathbb{T}_\mathbb{C}^\ell$ on a~complex vector space
$W$, it is well-known that the GIT-quotient $W/\!\!/\mathbb{T}_\mathbb{C}^\ell$ is isomorphic as a~complex algebraic
variety to a~complex orbifold if and only if the representation is \emph{simplicial} (see Subsection~\ref{reptheory}
and Section~\ref{simplex}).
By the Kempf--Ness theorem (to be recalled in Subsection~\ref{reptheory}), the symplectic quotient $M_0$ is
homeomorphic to such a~GIT-quotient.
Hence, the question arises whether the orbifold criterion in the complex algebraic setting carries over via the
Kempf--Ness homeomorphism to the symplectic setting.

Our results can be stated as follows.
If the symplectic quotient of a~unitary representation of a~compact torus is homeomorphic to an orbifold, then the
representation has to be simplicial (see Subsection~\ref{toricgeometry}).
We indicate methods of determining whether a~representation satisf\/ies this property directly from the weight matrix
in Section~\ref{simplex}.
This in particular resolves the conjectures stated in~\cite{HerbigIyengarPflaum}.
When the symplectic quotient has real dimension two, the representation is always simplicial; in this case, we further
demonstrate an explicit \emph{graded regular symplectomorphism} (to be def\/ined in Subsection~\ref{liftingthm}) to
a~quotient of~$\mathbb{C}$ by a~f\/inite abelian group.
On the other hand, we present in Subsection~\ref{finalargument} examples of simplicial unitary circle representations
whose symplectic quotients are homeomorphic to $\mathbb C^2$, for which there cannot exist a~graded regular
symplectomorphism to a~quotient of $\mathbb R^4$ by a~f\/inite subgroup of the group $\operatorname{Sp}(\mathbb R^4)$
of linear symplectomorphisms of $\mathbb R^4=T^*\mathbb R^2$.
So, roughly speaking, the simplicialness of the representation turns out to be merely a~necessary condition for the
existence of a~graded regular symplectomorphism with a~quotient by a~f\/inite group.

\looseness=-1
The reader might have noticed that our results should be taken with a~grain of salt.
Namely, for our counterexamples we cannot disprove the existence of a~symplectomorphism (see
Def\/i\-ni\-tion~\ref{symplectomorphism}) using the methods presented here, as the invariants we compute to distinguish
them from quotients by f\/inite groups are merely invariant under graded regular symplectomorphism.
More precisely, what we actually do is to focus on the case of real dimension 4 and work through the list of f\/inite
subgroups of the unitary group $\operatorname U_2$.
The Hilbert series of the ring of real polynomial invariants of these f\/inite subgroups are in principle computable by
Molien's formula, and we argue that the Hilbert series of the ring of regular functions on our symplectic circle
quotients cannot occur in this list.
This method is admittedly brute force, but it has the potential to guide us to a~classif\/ication of unitary symplectic
circle representations whose symplectic quotients are graded regularly symplectomorphic to quotients of unitary
representations of f\/inite groups.
We aim to complete this classif\/ication in the near future.
In higher dimensions, a~more intelligent approach is necessary.

Regular and graded regular symplectomorphism of singular phase spaces are roughly speaking those that can be obtained
using complete sets of dif\/ferentiable invariants.
In all practical applications, these are provided by the theorem of Schwarz--Mather~\cite{MatherDifferentiableInvariants,SchwarzSmoothFunctions} (see Theorem~\ref{diffinv}).
Though this construction principle for symplectomorphisms might look familiar to the specialist, we propose the
terminology in Section~\ref{smooth} to provide a~clear way of thinking about maps between singular phase spaces.
We expect that this language will have applications elsewhere.

\vspace{-2mm}

\section{Basic setup}
%\label{basics}
\vspace{-1mm}

\subsection{Background from representation theory}
\label{reptheory}

Here we recall some well-known facts about quotients of linear actions of compact groups and their relationship to
certain GIT-quotients.
For a~more systematic presentation we refer to G.W.~Schwarz' article~\cite{SchwarzTopAlgQuot}.

Let $G\to\operatorname{Gl}(W)$ be a~representation of a~compact Lie group on a~f\/inite-dimensional
real vector space~$W$.
By a~theorem of Hilbert and Hurwitz, there is a~complete system of real homogeneous polynomial invariants
$\rho_1,\dots,\rho_k$ in $\mathbb R[W]^G$; one can assume that the system is minimal.
This system, which we will refer to as a~\emph{Hilbert basis}, gives rise to a~map
$\rho=(\rho_1,\dots,\rho_k)\colon  W\to\mathbb R^k$, the corresponding \emph{Hilbert map}.
It is known that $\rho$ is proper and separates $G$-orbits.
The induced map $\cc\rho\colon  W/G\to\mathbb R^k$ will be referred to as the \emph{Hilbert embedding}.
By the Tarski--Seidenberg principle $X:=\operatorname{im}(\rho)\subset\mathbb R^k$ is a~semialgebraic set.
The gradients of the~$\rho_i$ can be used to calculate the inequalities that determine~$X$ (cf.~\cite[\S~6]{SchwarzTopAlgQuot}).
The Zariski closure~$\cc X$ of~$X$ is determined by the polynomial relations among the~$\rho_i$'s.
By def\/inition, a~function~$f$ on~$X$ is \emph{smooth} if it is the restriction $f=F_{|X}$ to $X$ of a~smooth function
$F\in\mathcal C^\infty(\mathbb R^k)$.
The algebra~$\mathcal C^\infty(X)$ of smooth functions on~$X$ is a~nuclear Fr{\'e}chet algebra (see,
e.g.,~\cite{PflaumBook}).

A key result for the analytic study of such an orbit space $W/G$ is the theorem of Schwarz and
Mather~\cite{MatherDifferentiableInvariants, SchwarzSmoothFunctions} on dif\/ferentiable invariants.

\begin{theorem}[G.W.~Schwarz, J.~Mather]
\label{diffinv}
With the notation above the pullback $\rho^*\colon \mathcal C^\infty(X)\to\mathcal C^\infty(W)^G$, $f\mapsto
f\circ\rho$ is split surjective onto the Fr{\'e}chet algebra $\mathcal C^\infty(W)^G$ of smooth invariants on $W$.
\end{theorem}

In~\cite{MatherDifferentiableInvariants, SchwarzSmoothFunctions}, the authors use Theorem~\ref{diffinv} to prove the
existence of a~complete set of dif\/ferentiable invariants for a~$G$-manifold using Mostov's embedding theorem, i.e.\
a generating set for the algebra of smooth $G$-invariant functions.
In the case of a~$G$-representation, a~complete set of dif\/ferentiable invariants is given by a~Hilbert basis.
Using the language of Section~\ref{smooth}, this theorem implies that the Hilbert embedding $\cc\rho$ is actually
a~dif\/feomorphism from the dif\/ferential space $(W/G,\mathcal C^\infty(W)^G)$ onto the dif\/ferential space
$(X,\mathcal C^\infty(X))$.

Now suppose $G\to\operatorname{U}(V)$ is a~unitary representation of the compact Lie group $G$ on a~f\/inite-dimensional complex vector space $V$ with hermitian scalar product $\langle\:,\:\rangle$.
By convention, $\langle\:,\:\rangle$ is complex antilinear in the f\/irst argument.
Note that we can make any symplectic representation of~$G$ unitary by using an invariant compatible complex structure.
In order to express equation~\eqref{eq:MomentMap} transparently, it will be convenient to express real polynomials
using complex coordinates.
Let $\cc V$ be the complex conjugate vector space of $V$, and then the identity map on $V$ induces a~complex antilinear
map ${}^{-}\colon  V\to\cc V$, $v\mapsto\cc v$.
The complex conjugation ${}^{-}$ extends to a~real structure on the algebra $\mathbb C[V\times\cc V]$, and the ring of
\emph{real regular functions on} $V$ is def\/ined to be the subring of invariants with respect to ${}^{-}$, i.e.\
$\mathbb R[V]:=\mathbb C[V\times\cc V]^-$.
It is of course isomorphic to the ring of regular functions on the real vector space $V_{\mathbb R}$ underlying~$V$.

The group $G$ acts on $\cc V$ by $\cc v\mapsto(g^{-1})^t\cc v$.
Letting $G$ act on $V\times\cc V$ diagonally, and observing that this action commutes with ${}^{-}$, we obtain an
action of $G$ on $\mathbb R[V]$ by $\mathbb R$-algebra automorphisms.
This action can be seen as coming from the obvious $\mathbb R$-linear $G$-action on $V_\mathbb R$.
Hence $\mathbb R[V]^G$ is a~$\mathbb Z$-graded Noetherian $\mathbb R$-algebra, we can f\/ind a~Hilbert basis
$\rho_1,\dots,\rho_k\in\mathbb R[V]^G$ and Theorem~\ref{diffinv} applies.
Note that $v\mapsto\langle v,v\rangle$ is always a~quadratic invariant.

It is well-known that the unitary action of $G$ on $V$ extends uniquely to a~$\mathbb C$-linear action of the
complexif\/ication $G_{\mathbb C}$ of $G$ on $V$.
Note also that the complexif\/ication of the $G$-action on $V\times\cc V$ turns out to be the cotangent lifted
$G_{\mathbb C}$-action on $V\times V^*$.
Moreover, we have the following isomorphism of invariant rings
\begin{gather}
\label{invcotangent}
\mathbb R[V]^G\otimes_{\mathbb R}\mathbb C\cong\mathbb C[V\times V^*]^{G_{\mathbb C}}
\end{gather}
as $\mathbb Z$-graded $\mathbb C$-algebras.

The (inf\/initesimal) information of the unitary representations $G\to\operatorname U(V)$ can be encoded into the
moment map $J$.
This is the regular quadratic map
\begin{gather}
\label{eq:MomentMap}
J\colon \  V\to\mathfrak g^*,
\qquad
J_\xi(v)=\langle J(v),\xi\rangle:=\frac{\sqrt{-1}}{2}\langle v,\xi v\rangle
\end{gather}
for $\xi\in\mathfrak g$.
Alternatively, we can think of $J$ as a~linear map $\mathfrak g\to\mathbb R[V]$.
Often it is convenient to write $J_i:=J_{e_i}$ for some f\/ixes basis $e_1,\ldots,e_{\ell}$ of $\mathfrak g$.
The moment map is of particular importance when it comes to discussing the symplectic geometry of our unitary
representation.
Let us, for convenience, identify $V$ with $\mathbb C^n$ by choosing an orthonormal basis and denote the corresponding
coordinates by $(\boldsymbol z,\cc{\boldsymbol z})=(z_1,\dots,z_n,\cc z_1,\dots,\cc z_n)$.
It follows that $\mathbb R[V]$ is identif\/ied with $\mathbb R[\mathbb C^n]=\mathbb C[\boldsymbol z,\cc{\boldsymbol
z}]^-$.
The Poisson bracket corresponding to the symplectic form $\omega\in\Omega^2(V)$, $\omega(v,w)=\operatorname{Im}\langle
v,w\rangle$, is given by the relation
\begin{gather*}
\{z_i,\cc z_j\}=\frac{2}{\sqrt{-1}}\delta_{i,j},
\end{gather*}
all other brackets between coordinates being zero.
This makes $\mathcal C^\infty(\mathbb C^n)$ into a~Poisson algebra with Poisson subalgebra $\mathbb R[\mathbb C^n]$.
It turns out that $\{J_\xi,J_\eta\}=J_{[\xi,\eta]}$, which is equivalent to the equivariance of the map
$J\colon  V\to\mathfrak g^*$.

In the situation of a~unitary representation the zero f\/ibre $Z=J^{-1}(0)$ of the moment map always has a~conical
singularity at $0$.
The symplectic quotient $M_0=Z/G$ is a~stratif\/ied symplectic space (this will be further explained in
Section~\ref{smooth}).
In general, it is not a~real variety but a~semialgebraic set.
In contrast, the GIT quotient $V/\!\!/G_{\mathbb C}$ is def\/ined to be the complex variety underlying the $\mathbb
C$-algebra $\mathbb C[V]^{G_{\mathbb C}}$.
It might happen that $V/\!\!/G_{\mathbb C}$ is actually smooth (cf.\
Section~\ref{counter}).
Due to the following theorem of Kempf and Ness (see~\cite[Corollary~4.7]{SchwarzTopAlgQuot}), $Z=J^{-1}(0)$ is
sometimes called the \emph{Kempf--Ness set}.

\begin{theorem}[G.~Kempf and L.~Ness]
The map $Z\hookrightarrow V\mapsto V/\!\!/G_{\mathbb C}$ is proper and induces a~homeomorphism $Z/G\to
V/\!\!/G_{\mathbb C}$.
\end{theorem}

In view of equation~\eqref{invcotangent}, the Kempf--Ness theorem actually comes as a~surprise, as the invariant theory
of a~cotangent lifted representation is more involved than that of the original representation.
The theorem is a~useful tool to count dimensions of symplectic quotients.
The aim of the paper is to give examples where $V/\!\!/G_{\mathbb C}$ is smooth while $Z/G$ is not an orbifold in an
appropriate sense.

\subsection{Background from toric geometry}
\label{toricgeometry}

Next we would like to specialize the discussion to the case when our compact group $G$ is actually an
$\ell$-dimensional torus.
By this we mean an $\ell$-fold copy $\mathbb{T}^\ell:=(\mathbb{S}^1)^\ell$ of the unit sphere
$\mathbb{S}^1\subset\mathbb C$.
We are interested in unitary representations
\begin{gather*}
G=\mathbb{T}^\ell\to\operatorname U_n:=\operatorname U(\mathbb C^n),
\end{gather*}
where $\mathbb C^n$ is understood with its standard hermitian scalar product as in the previous section.
We identify the Lie algebra $\mathfrak g$ of $G=\mathbb{T}^\ell$ with $\mathbb R^\ell$ by writing an arbitrary element
$(t_1,\dots,t_\ell)\in G=\mathbb{T}^\ell$ in the form $t_i=\exp(2\pi\sqrt{-1}\xi_i)$, for the vector
$(\xi_1,\ldots,\xi_\ell)\in\mathfrak g=\mathbb R^n$.
Since the factors $\mathbb{S}^1$ of our torus action can be simultaneously diagonalized, the unitary representation can
actually be encoded into a~\emph{weight matrix} $A=(a_{ij})\in\mathbb{Z}^{\ell\times n}$.
More specif\/ically, setting $(\eta_1,\ldots,\eta_n):=(\xi_1,\ldots,\xi_\ell)\cdot A\in\mathbb R^n$, the
$G=\mathbb{T}^\ell$-action corresponding to the weight matrix $A$ is given by the formula
\begin{gather*}%\label{eq-IntroLieAlgAction}
(t_1,\ldots,t_\ell).
(z_1,\ldots,z_n)=\big(\exp\big(2\pi\sqrt{-1}\eta_1\big)z_1,\ldots,\exp\big(2\pi\sqrt{-1}\eta_n\big)z_n\big).
\end{gather*}
Elementary row operations with integer scalars for $A$, i.e.\
row operations that correspond to left multiplication by elements of $\mathrm{GL}_\ell(\mathbb{Z})$, correspond to the
changing of a~basis of $\mathfrak{g}$, while permutations of the columns of $A$ correspond to changing coordinates for
$\mathbb{C}^n$.

The components $J_i:=J_{e_i}=\langle J,e_i\rangle$ of the moment map
$J\colon \mathbb{C}^n\to\mathbb{R}^\ell\cong\mathfrak g^\ast$ can also be expressed in terms of the weight
matrix
\begin{gather*}%\label{eq:MomentMapComps}
J_i(\boldsymbol z,\cc{\boldsymbol z})=\frac{1}{2}\sum\limits_{j=1}^n a_{ij}z_j\cc z_j,
\qquad
i=1,\ldots,\ell.
\end{gather*}
\begin{note}
\label{not:zAxAyA}
Sometimes it will be convenient to emphasize the dependency on $A$ in the notation.
In these cases we will write $J_A$ for the moment map, $Z_A=J_A^{-1}(0)$ for the zero f\/iber, and
$M_A=Z_A/\mathbb{T}^\ell$ for the reduced space.
We will also let $X_A=Z_A\cap\mathbb{S}^{2n-1}$ denote the intersection of the zero f\/iber with the unit sphere in
$\mathbb{C}^n$ and $Y_A=X_A/\mathbb{T}^\ell$ the link.
Note that $X_A$ is clearly $\mathbb{T}^\ell$-invariant.
\end{note}

The case of toric moment maps has certain peculiarities; for example, the components of toric moment maps are actually
invariants.
We will have to say more about this in Subsection~\ref{magic}.

Let us introduce some further notation.
We denote by $\operatorname{sq}\colon \mathbb{C}^n\to\mathbb{R}^n$ the map $(z_1,\ldots,z_n)\mapsto(z_1\cc
z_1,\ldots,z_n\cc z_n)$.
It is clear that $\operatorname{sq}$ is actually $G=\mathbb{T}^\ell$-invariant and hence induces a~map
$\widetilde{\operatorname{sq}}\colon \mathbb{C}^n/\mathbb{T}^\ell\to\mathbb{R}^n$.

We will primarily be interested in the case where the action of $\mathbb{T}^\ell$ on $\mathbb{C}^n$ is ef\/fective, i.e.\
if for some $t\in\mathbb{T}^\ell$ we have $tz=z$ for all $z\in\mathbb{C}^n$, then $t=1$.
We will see below (Lemma~\ref{effectivemoment}) that this introduces no loss of generality.
In order to do so, we f\/irst interpret this condition in terms of the weight matrix $A$.

It is easy to see that there is a~subgroup $K\leq\mathbb{T}^\ell$ of positive dimension that acts trivially on~$\mathbb{C}^n$ if and only if $\operatorname{rank}(A)<\ell$.
In particular, choosing a~basis for $\mathfrak{g}$ that contains an element of the Lie algebra $\mathfrak{k}$ of~$K$,
it is easy to see that the corresponding row of~$A$ is the zero row.
Hence, $A$ has full rank if and only if the subgroup of~$\mathbb{T}^\ell$ that acts trivially on~$\mathbb{C}^n$ is
f\/inite.
In this case, we have the following lemma; we include the proof since we do not know of an appropriate reference.

\begin{lemma}
%\label{effectivemat}
Suppose $A\in\mathbb{Z}^{\ell\times n}$ has full rank $\ell\leq n$.
Then the action of $\mathbb{T}^\ell$ on $\mathbb{C}^n$ is effective if and only if the nonzero
$\ell\times\ell$-minors of $A$ are relatively prime.
Moreover, if~$p$ is a~prime that divides each of the $\ell\times\ell$-minors of~$A$, by elimination with integer
scalars and permuting columns, $A$~can be expressed in a~form where each entry of its first row is divisible by~$p$.
\end{lemma}

\begin{proof}
Suppose $t=\big(\exp(2\pi\sqrt{-1}\xi_1),\ldots,\exp(2\pi\sqrt{-1}\xi_\ell)\big)\in\mathbb{T}^\ell$ is nontrivial and acts
trivially on $\mathbb{C}^n$.
As $A$ has full rank, $t$ must be of f\/inite order.
Thus there is a~$j\in\{1,\dots,\ell\}$ such that $\xi_j=k/q$ for some coprime integers $k,q$ with $q\ge2$.
By assumption, $(\xi_1,\ldots,\xi_\ell)A\in\mathbb Z^n$.
Let $B$ be a~nonsingular $\ell\times\ell$-submatrix of $A$.
Since $(\xi_1,\ldots,\xi_\ell)B\in\mathbb Z^\ell$, we can use Cramer's rule to conclude that $q\mid\det(B)$.

Conversely, let $\gcd_\ell(A)$ denote the $\gcd$ of the $\ell\times\ell$-minors of $A$.
We prove by induction on $n-\ell$ that if $p$ is a~prime that divides $\gcd_\ell(A)$, then $A$ can be row-reduced with
integer scalars and the coordinates of $\mathbb{C}^n$ can be permuted so that $(\exp(2\pi\sqrt{-1}/p),0,\ldots,0)$ acts
trivially on~$\mathbb{C}^n$, i.e.\
$(1/p,0,\ldots,0)A\in\mathbb Z^n$.

Let $A\in\mathbb Z^{\ell\times n}$ and let $p$ be a~prime such that $p|\gcd_\ell(A)$.
Assume the result holds for all $\ell\times(\ell+k)$-weight matrices with $k<n-\ell$.
By row operations and permutations of coordinates, we can assume that $A=[D\:|\:C]$ where
$D=\operatorname{diag}(d_1,\ldots,d_\ell)\in\mathbb Z^{\ell\times\ell}$ and $C\in\mathbb{Z}^{\ell\times(n-\ell)}$.
Then $p|\det(D)$ so that by further permuting coordinates, we can assume that $p|d_1$.
If $n-\ell=0$ it follows that $(1/p,0,\ldots,0)A\in\mathbb{Z}^n$.

Otherwise, let $A^\prime$ denote the matrix formed by removing the f\/irst column of $A$.
Consider the case when $A^\prime$ does not have full rank.
This means in particular that each $\ell\times\ell$-submatrix of $A^\prime$ corresponding to the columns
$2,3,\ldots,\ell,\ell+j$ of $A$ for $1\leq j\leq n-\ell$ is singular.
It follows that the f\/irst row of $A^\prime$ is the zero row, which implies $(1/p,0,\ldots,0)A\in\mathbb{Z}^n$.

On the other hand, suppose $A^\prime$ has full rank.
Then by the inductive hypothesis, we can row-reduce $A^\prime$ with integer scalars and permute the columns of
$A^\prime$ to yield a~matrix $R^\prime$ such that that $(1/p,0,\ldots,0)R^\prime\in\mathbb{Z}^{n-1}$.
If we apply the same row-reduction to $A$, however, and permute columns $2,3,\ldots,n$ in the same way, it is easy to
see that the resulting matrix $R$ is of the form $[\boldsymbol{r}\:|\:R^\prime]$ where
$R^\prime\in\mathbb{Z}^{\ell\times(n-1)}$ and $\boldsymbol{r}$ is a~column with each entry divisible by $p$.
It follows that $(1/p,0,\ldots,0)R\in\mathbb{Z}^n$, completing the proof.
\end{proof}

Now, suppose the action of $\mathbb{T}^\ell$ on $\mathbb{C}^n$ is not ef\/fective, and let $K\leq\mathbb{T}^\ell$
denote the subgroup that acts trivially.
Then $\mathbb{T}^\ell$ f\/ibers over $\mathbb{T}^\ell/K$, which is itself a~torus, and we may consider the Hamiltonian
action of $\mathbb{T}^\ell/K$ on $\mathbb{C}^n$.
If $K$ is inf\/inite and connected, then $\mathbb{T}^\ell/K$ is a~torus of dimension smaller than $\ell$.
The row-reduced weight matrix $A$ has zero rows, and the moment maps of the $\mathbb{T}^\ell$- and
$\mathbb{T}^\ell/K$-actions dif\/fer only by extending by zero.
If $K$ is f\/inite, then the moment maps of the two actions coincide up to an isomorphism between the Lie algebra of~$\mathbb{T}^\ell/K$ with that of $\mathbb{T}^\ell$.
Combining these two arguments for an arbitrary~$K$ yields the following.

\begin{lemma}
\label{effectivemoment}
Let $A^\prime$ denote the weight matrix of the $\mathbb{T}^\ell/K$-action on~$\mathbb{C}^n$.
Then $J_A^{-1}(0)=J_{A^\prime}^{-1}(0)$.
\end{lemma}

As a~consequence, if the action of $\mathbb{T}^\ell$ is not ef\/fective, then we may replace $\mathbb{T}^\ell$ with
$\mathbb{T}^\ell/K$ without changing the reduced space.
Hence, in the sequel, we assume without loss of generality that $\mathbb{T}^\ell$ acts ef\/fectively on $\mathbb{C}^n$,
and in particular that $\ell\leq n$.

Now, let $\mathbb{T}_\mathbb{C}^\ell$ denote the complexif\/ication of $\mathbb{T}^\ell$.
Then $\mathbb{T}_\mathbb{C}^\ell$ acts on $\mathbb{C}^n$ via
\begin{gather*}
(w_1,\ldots,w_\ell)(z_1,\ldots,z_n)=
\big(w_1^{a_{11}}w_2^{a_{21}}\cdots w_\ell^{a_{\ell1}}z_1,\ldots,w_1^{a_{1n}}w_2^{a_{2n}}\cdots w_\ell^{a_{\ell n}}z_n\big),
\end{gather*}
and this action induces an injective homomorphism $\mathbb{T}_\mathbb{C}^\ell\to\mathbb{T}_\mathbb{C}^n$.
Then the GIT quotient $\mathbb{C}^n/\!\!/\mathbb{T}_\mathbb{C}^\ell$ is equipped with an ef\/fective action of
$\mathbb{T}_\mathbb{C}^n/\mathbb{T}_\mathbb{C}^\ell\cong\mathbb{T}_\mathbb{C}^{n-\ell}$ with a~single, dense orbit and
hence has the structure of an $(n-\ell)$-dimensional toric variety $\mathcal{X}$, see e.g.~\cite{CoxLittleSchenck} or~\cite{FultonBook}.
In particular, $\mathbb{C}^n/\!\!/\mathbb{T}_\mathbb{C}^\ell$ is the af\/f\/ine toric variety given by the spectrum of
the semigroup $\ker(A)\cap\mathbb{Z}_{\geq0}^n$ and hence is associated to the cone given by the kernel of $A$
intersected with the positive $n$-ant in~$\mathbb{R}^n$.
\begin{definition}
\label{def:SigmaA}
The cone $\sigma_A$ associated to the weight matrix $A$ is the intersection of the kernel of $A$ with the positive
$n$-ant in $\mathbb{R}^n$.
\end{definition}

Recall that the cone $\sigma_A$ is \emph{simplicial} if it is generated by a~collection of linearly independent vectors.
It is well-known, see e.g.~\cite[Section~2.2]{FultonBook}, that if $\sigma_A$ is simplicial, then the af\/f\/ine toric variety associated to
$\sigma_A$ is a~complex orbifold.
In particular, applying the Cox construction, see~\cite[Chapter~5]{CoxLittleSchenck}, we have that that
$\mathcal{X}=\mathbb{C}^{n-\ell}/\Gamma$ for a~f\/inite group $\Gamma$ as follows.

We have the short exact sequence~\cite[Theorem~4.1.3]{CoxLittleSchenck}
\begin{gather*}
0\longrightarrow M\longrightarrow
\operatorname{Div}_{\mathbb{T}_\mathbb{C}^{n-\ell}}(\mathcal{X})\longrightarrow
\operatorname{Cl}(\mathcal{X})\longrightarrow0,
\end{gather*}
where $M$ denotes the character lattice of the algebraic torus $\mathbb{T}_\mathbb{C}^{n-\ell}$,
$\operatorname{Div}_{\mathbb{T}_\mathbb{C}^{n-\ell}}(\mathcal{X})$ denotes the group of
$\mathbb{T}_\mathbb{C}^{n-\ell}$-invariant Weil divisors of $\mathcal{X}$, and $\operatorname{Cl}(\mathcal{X})$ denotes
the class group of~$\mathcal{X}$.
Choosing bases, this sequence can be expressed as
\begin{gather*}
0\longrightarrow\mathbb{Z}^{n-\ell}\stackrel{(\ast)}{\longrightarrow}\mathbb{Z}^{n-\ell}\longrightarrow
\operatorname{Cl}(\mathcal{X})\longrightarrow0,
\end{gather*}
where the map $(\ast)$ is given by the matrix whose rows are the coordinates of the $n-\ell$ linearly independent
minimal generators of the cone $\sigma_A$ and hence has maximal rank.
In particular, $\operatorname{Cl}(\mathcal{X})$ is f\/inite.
Applying $\operatorname{Hom}_\mathbb{Z}(\cdot,\mathbb{T}_\mathbb{C}^1)$ and setting
$\Gamma=\operatorname{Hom}_\mathbb{Z}(\operatorname{Cl}(\mathcal{X}),\mathbb{T}_\mathbb{C}^1)$ yields the exact sequence
\begin{gather*}
1\longrightarrow\Gamma\longrightarrow\mathbb{T}_\mathbb{C}^{n-\ell}
\stackrel{(\ast)^T}{\longrightarrow}\mathbb{T}_\mathbb{C}^{n-\ell}\longrightarrow1,
\end{gather*}
def\/ining an action of $\Gamma$ on $\mathbb{C}^{n-\ell}$.
Hence, as $\sigma_A$ consists of a~single cone so that the exceptional set is empty, the toric variety $\mathcal{X}$ is
given by the complex orbifold $\mathbb{C}^{n-\ell}/\Gamma$.

In particular, if $n-\ell=1$, it is easy to see that the cone $\sigma_A$ is simply $\mathbb{R}_{\geq0}$ with minimal
generator $1$.
Therefore, the map $\mathbb{Z}\stackrel{(\ast)}{\longrightarrow}\mathbb{Z}$ above is simply the identity, and
$\operatorname{Cl}(\mathcal{X})$ and $\Gamma$ are both trivial.
It follows that $\mathcal{X}=\mathbb{C}$.

For any complex orbifold $Q$, each local group action preserves the complex structure.
It follows that $Q$ is a~\emph{locally orientable} orbifold, i.e.\
each local group action preserves a~local orientation.
By~\cite[4.2.4]{Haefliger}, the underlying topological space of a~locally orientable orbifold of (real) dimension~$m$
is an \emph{$m$-dimensional rational homology manifold}.
That is, if $\mathbb{X}_Q$ denotes the underlying space of $Q$, then the local homology groups with rational
coef\/f\/icients $H_k(\mathbb{X}_Q,\mathbb{X}_Q-x;\mathbb{Q})$ at each point $x\in\mathbb{X}_Q$ satisfy
\begin{gather*}
H_k(\mathbb{X}_Q,\mathbb{X}_Q-x;\mathbb{Q})=
\begin{cases} \mathbb{Q}, &  k=m,\\
0,& k\neq m.
\end{cases}
\end{gather*}

If, on the other hand, $\sigma_A$ is not simplicial, then the recursion formula given in~\cite[p.~2]{BarBraFieKau}
for the local intersection cohomology Betti numbers in terms of the cone generators indicates that the second local
intersection cohomology is nontrivial.
Because the local intersection cohomology of a~rational homology manifold is trivial, it follows that the toric variety~$\mathcal{X}$ associated to $\sigma_A$ is not a~rational homology manifold.
With this, applying the Kempf--Ness homeomorphism between $M_A=Z_A/\mathbb{T}^\ell$ and $\mathcal{X}$, we have the
following.
\begin{theorem}
\label{orbifoldsimplicial}
Using the notation of Note~{\rm \ref{not:zAxAyA}} and Definition~{\rm \ref{def:SigmaA}}, the reduced space
$M_A=Z_A/\mathbb{T}^\ell$ associated to $A\in\mathbb{Z}^{\ell\times n}$ is a~rational homology manifold if and only if
the cone $\sigma_A$ is simplicial.
\end{theorem}

In particular, note that symplectic orbifolds are locally orientable and hence rational homo\-lo\-gy manifolds.
Therefore, if the cone~$\sigma_A$ is not simplicial, then the topological space~$M_A$ does not admit a~homeomorphism to
a~symplectic orbifold.

In the sequel, it will be convenient to use the following terminology.
\begin{definition}\label{def:Simplicial}
We say that the weight matrix $A\in\mathbb Z^{\ell\times n}$ is \emph{simplicial} if the corresponding cone~$\sigma_A$
is simplicial.
In this case we also say that the corresponding unitary $\mathbb{T}^\ell$-action and its complexif\/ied
$\mathbb{T}^\ell_{\mathbb C}$-action are \emph{simplicial}.
\end{definition}

\subsection{Other topological indications}
%\label{topology2}

In many examples of non-simplicial weight matrices $A$, it is possible to demonstrate that the reduced space is not
homeomorphic to a~symplectic orbifold directly without appealing to the Kempf--Ness homeomorphism.
In this subsection, we brief\/ly indicate results in this direction.

In~\cite[Example~2.4]{CushmanSjamaar}, the reduced space corresponding to the weight matrix $[-1,-1,1,1]$ was described
as the cone on $\mathbb{S}^3\times_{\mathbb{S}^1}\mathbb{S}^3$, implying that the local homology in degree $3$ at the
cone point is nontrivial.
It follows that the reduced space is not a~rational homology manifold and hence not an orbifold.

By~\cite[Proposition~3.1]{HatcherVogtmann}, the quotient of an $n$-dimensional sphere by a~f\/inite group acting
linearly and preserving orientation is a~rational homology $n$-sphere, i.e.\
has the homology with rational coef\/f\/icients of the $n$-dimensional sphere $\mathbb{S}^n$.
It follows that the link $Y_A=X_A/\mathbb{T}^\ell$, see Note~\ref{not:zAxAyA}, of a~locally orientable $n$-dimensional
orbifold singularity is a~rational homology $n$-sphere.
In~\cite{HerbigIyengarPflaum}, this observation was used to show that the reduced space $M_A=Z_A/\mathbb{T}^\ell$
cannot be an orbifold if the link $Y_A$ is not a~rational homology sphere.
In particular, in the case $\ell=1$,~\cite[Proposition~3.1]{HerbigIyengarPflaum} demonstrates that $Y_A$ is not
a~rational homology sphere if the weight matrix $A$ has at least two positive and two negative entries; this condition
is clearly equivalent to the negation of Theorem~\ref{simpcondgen}(2) below in this case.

Similarly, in cases where $X_A$ consist of points of a~single orbit type, the quotient map~$X_A\to Y_A$ is a~torus
f\/ibration with f\/iber given by the quotient of $\mathbb{T}^\ell$ by the isotropy group of~$X_A$.
In this case, formulas for the homology of~$X_A$ have been developed in~\cite{BosioMeersseman}, and in some cases, the
exact sequence~\cite[Theorem~2, p.~482]{SpanierBook} can be used to demonstrate that~$X_A$ does not admit such
a~torus f\/ibration over a~rational homology sphere of the appropriate dimension.

More generally, note that this argument can be applied to the closed orbit-type strata of the link $Y_A$ to show that
the reduced space does not admit a~stratum-preserving homeomorphism to an orbifold.
To see this, suppose $G$ is a~f\/inite group acting on a~sphere~$\mathbb{S}^n$.
For each $H\leq G$, we let $\mathbb{S}^n_H$ denote the set of points with isotropy group $H$ and $\mathbb{S}^n_{(H)}$
the set of points with orbit type~$(H)$.
Then $H$ acts trivially on $\mathbb{S}^n_H$, $N_G(H)/H$ acts freely on $\mathbb{S}^n_H$, and $\mathbb{S}^n_{(H)}/G$ is
dif\/feomorphic to $\mathbb{S}^n_H/N_G(H)$; see~\cite[Theorem~4.3.10 and Corollary~4.3.11]{PflaumBook}.
If $\mathbb{S}^n_{(H)}$ has minimal dimension among the strata, then $\mathbb{S}^n_H=(\mathbb{S}^n)^H$, and hence
$\mathbb{S}^n_{(H)}/G$ is dif\/feomorphic to the quotient of a~sphere by the free action of a~f\/inite group.
Similarly, if $\mathbb{S}^n_H$ is closed, then it is locally a~stratum of minimal dimension, and we can draw the same
conclusion.
It follows that the closed orbit-type strata of the link of an orbifold singularity are as well rational homology
spheres, so that this must also be true for a~reduced space that admits a~stratum-preserving homeomorphism to an
orbifold.

We illustrate these observations with the following.
\begin{example}
\label{nonsimpexample}
Consider the case of $\mathbb{T}^2$ acting on $\mathbb{C}^6$ with weight matrix
\begin{gather*}
A=\left[
\begin{matrix}
1&-1&1&-1&0&0
\\
0&0&0&0&1&-1
\end{matrix}
\right].
\end{gather*}
Then $Z_A$ is described by
\begin{gather*}
|z_1|^2+|z_3|^2=|z_2|^2+|z_4|^2,
\qquad
|z_5|^2=|z_6|^2.
\end{gather*}
The isotropy types away from the origin are given by $(z_1,z_2,z_3,z_4,0,0)$ with isotropy $1\times\mathbb{T}^1$,
$(0,\ldots,0,z_5,z_6)$ with isotropy $\mathbb{T}^1\times1$, and $(z_1,\ldots,z_6)$ with trivial isotropy.
If $z_5=z_6=0$, then the intersection with the unit sphere is $|z_1|^2+|z_2|^2=|z_3|^2+|z_4|^2=1/2$, and the
corresponding orbit-type stratum is homeomorphic to $\mathbb{S}^3\times\mathbb{S}^3/\mathbb{T}^1$.
Using~\cite[Theorem~2, p.~482]{SpanierBook}, it is an easy exercise to show that $\mathbb{S}^3\times\mathbb{S}^3$
does not admit a~$\mathbb{T}^1$-f\/ibration over a~rational homology $5$-sphere, and hence that a~closed stratum of
$Y_A=X_A/\mathbb{T}^2$ is not a~rational homology $5$-sphere.
It follows that $M_A=Z_A/\mathbb{T}^2$ does not admit a~stratum-preserving homeomorphism with an orbifold.
\end{example}

\section{Gaussian elimination and the simplicial condition}
\label{simplex}

In this section, we will use Theorem~\ref{orbifoldsimplicial} to determine necessary and suf\/f\/icient conditions for
the reduced space $M_A=Z_A/\mathbb{T}^\ell$ to be a~rational homology manifold directly in terms of the matrix~$A$.
Given a~subset $X$ of $\mathbb R^n$ we write $\operatorname{af\/f}(X)$ for its af\/f\/ine hull and
$\operatorname{cch}(X)$ for its closed convex hull.
By~$X^\circ$ we mean its relative interior, i.e., the interior of~$X$ seen as a~subspace of $\operatorname{af\/f}(X)$.
We also use the shorthand~$X^c$ for the complement $\mathbb R^n\backslash X$.

Let $A\in\mathbb{Z}^{\ell\times n}$ with $\ell\leq n$.
Let $\Delta^{n-1}$ denote the standard simplex in $\mathbb{R}^n$, and let
$\mathcal{P}_A:=\operatorname{ker}(A)\cap\Delta^{n-1}\subset\mathbb{R}^n$ denote the intersection of the kernel of~$A$
in $\mathbb{R}^n$ with the standard simplex~$\Delta^{n-1}$.
Then the cone $\sigma_A$ def\/ined in Def\/inition~\ref{def:SigmaA} is spanned by~$\mathcal{P}_A$.
Note that if $\mathcal{P}_A\neq\varnothing$, then $\mathcal{P}_A$ is a~polytope by~\cite[Corollary~9.4]{BronstedBook}.
Each element of~$\mathcal{P}_A$ is a~convex combination of its vertices by def\/inition, so that the vertices of~$\mathcal{P}_A$ clearly span the linear space spanned by~$\mathcal{P}_A$.
It follows that if~$\mathcal{P}_A$ has dimension~$m$, then the vertices are linearly independent if and only if there
are exactly $m+1$ vertices.
This is the case if and only if $\mathcal{P}_A$ is combinatorially equivalent to a~standard simplex, see~\cite[Chapter~2,
\S~10]{BronstedBook}, so that the matrix~$A$ is simplicial if and only if~$\mathcal{P}_A$ is combinatorially
equivalent to a~simplex.

In examples, the most direct method of determining whether $A$ is simplicial is to compute the vertices of
$\mathcal{P}_A$ using the results of Lemma~\ref{vertices} below.
However, in the sequel, we will need to use a~standard row-reduced form of a~simplicial weight matrix $A$, and hence we
will reformulate the simplicial condition (cf.
Def\/inition~\ref{def:Simplicial}) in these terms in Theorem~\ref{simpcondgen}.
In addition, we give a~geometric formulation to aid in the reader's intuition.

We use $e_1,\ldots,e_n$ to denote the standard basis vectors of $\mathbb{R}^n$ so that
$\Delta^{n-1}=\operatorname{cch}(e_1,\ldots,e_n)$ is the closed convex hull of the set of standard basis vectors.
If $I\subset\{1,\ldots,n\}$ is a~subset of indices, we let
\begin{gather*}
V_I=\big\{(x_1,\ldots,x_n)\in\mathbb{R}^n\mid x_j=0
\
\forall\, j\in\{1,\ldots,n\}\setminus I\big\}
\end{gather*}
denote the coordinate subspace associated to $I$.
Recall that $X_A=Z_A\cap\mathbb{S}^{2n-1}$ denotes the intersection of the zero f\/iber $Z_A$ with the unit sphere in
$\mathbb{C}^n$ and $Y_A=X_A/\mathbb{T}^\ell$ denotes the link, see Note~\ref{not:zAxAyA}.
Then we have that $\operatorname{sq}(X_A)=\widetilde{\operatorname{sq}}(Y_A)=\mathcal{P}_A$, where $\operatorname{sq}$
and $\widetilde{\operatorname{sq}}$ are the maps def\/ined in Subsection~\ref{toricgeometry}.
As well, note that the combinatorial type of $\mathcal{P}_A$ is clearly invariant under row reduction and permuting the
columns of $A$.

In general, it may happen that $\mathcal{P}_A$ is contained in a~coordinate subspace of $\mathbb{R}^n$ and hence
a~proper face of $\Delta^{n-1}$.
To address this possibility, let
\begin{gather*}
I_A=\big\{j\in\{1,\ldots,n\}\, |\,\exists\, (x_1,\ldots,x_n)\in\mathcal{P}_A:x_j\neq0\big\}
\end{gather*}
denote the set of coordinates $x_j$ that are not identically $0$ on $\mathcal{P}_A$.
Equivalently, $I_A$ is the set of indices $j$ such that there is an element of $\operatorname{ker}(A)$ with
non-negative entries and positive $j$th entry.
Let $A^\prime$ denote the $\ell\times|I_A|$ submatrix of $A$ given by the columns corresponding to elements of~$I_A$.
Let $V_{I_A}$ denote the coordinate subspace of $\mathbb{R}^n$ associated to~$I_A$, i.e., the subspace
$\{(x_1,\ldots,x_n)\in\mathbb{R}^n\mid x_j=0
\
\forall\, j\notin I_A\}$.
In examples, $I_A$ can be determined by computing the vertices of~$\mathcal{P}_A$.
Let $r\leq\ell$ denote the rank of~$A$.
We will establish the following criteria for the a~simplicial weight matrix.
Note that $\langle\cdot,\cdot\rangle$ denotes the standard inner product on $\mathbb{R}^n$.
\begin{theorem}\label{simpcondgen}
Let $A$ be an $n\times\ell$ weight matrix.
The following are equivalent.
\begin{enumerate}\itemsep=0pt
\item[$1.$] The polytope $\mathcal{P}_A$ is combinatorially equivalent to a~simplex.
\item[$2.$] By permuting the indices in $I_A$ and performing elementary row operations with integer scalar multiples, the
matrix $A^\prime$ can be expressed in the form
\begin{gather*}
\left[
\begin{matrix}
D
\\
0
\end{matrix}
\right|\left.
\begin{matrix}C&0
\\
0&0
\end{matrix}
\right],
\end{gather*}
where $D$ is an $r\times r$ diagonal matrix with strictly negative entries on the diagonal and $C$ is an $r\times q$
matrix such that each entry is nonnegative and $q\leq|I_A|-r$.
\item[$3.$] There are vectors $\mu_1,\ldots,\mu_r\in V_{I_A}$ and indices $j_1,\ldots,j_r\in I_A$ such that
$\operatorname{ker}(A^\prime)=V_{I_A}\cap\left(\bigcap_{i=1}^r\mu_i^\perp\right)$, where $\mu_i^\perp$ denotes the
orthogonal complement in $\mathbb{R}^n$, and for each $i=1,\ldots,r$, $\langle e_{j_i},\mu_i\rangle<0$, $\langle
e_{j_k},\mu_i\rangle=0$ for $k\neq i$, and $\langle e_j,\mu_i\rangle\geq0$ for $j\in I_A$, $j\neq j_i$.
\end{enumerate}
These conditions are trivially satisfied if $n\leq r+2$.
\end{theorem}

Note that in condition (2) of Theorem~\ref{simpcondgen}, by construction of the index set $I_A$, the matrix $C$ cannot
have rows that are identically zero.

Condition (3) of Theorem~\ref{simpcondgen} can be understood as follows.
For each $i=1,\ldots,r$, let $H_i=\mu_i^\perp\cap V_{I_A}$ denote the orthogonal complement $\mu_i^\perp$ in $V_{I_A}$.
Then condition (3) states that the hyperplane $H_i$ separates one vertex of the standard simplex in $V_{I_A}$ from the
others, and moreover that each hyperplane $H_i$ contains all of the separated basis vectors $e_{j_k}$ for $k\neq i$.
This condition lends some intuition for the geometric meaning of simplicial condition (2).

In order to establish Theorem~\ref{simpcondgen}, we will f\/irst restrict to the case of weight matrices satisfying the
following hypotheses to simplify the arguments.
\begin{enumerate}[$(i)$]\itemsep=0pt
\item The polytope $\mathcal{P}_A$ has nonempty intersection with the relative interior of the standard simplex
$\Delta^{n-1}$.
\item The matrix $A$ has full rank $\ell$.
\item The matrix $A$ has no columns that are identically zero.
\end{enumerate}
Note that as each point in $(\Delta^{n-1})^\circ$ has nonzero $x_i$-coordinate for each $i$, hypothesis $(i)$ is
equivalent to $I_A=\{1,\ldots,n\}$.
Similarly, $(ii)$~implies that $\operatorname{ker}(A)$ has dimension $n-\ell$; equivalently, no positive-dimensional
subgroups of $\mathbb{T}^\ell$ act trivially on $\mathbb{C}^n$.
Hypothesis $(iii)$~implies that there are no coordinate lines in $\mathbb{C}^n$ on which $\mathbb{T}^\ell$ acts trivially.
Assuming $(i)$, $(ii)$, and $(iii)$, it is easy to see that the relative interior of $\mathcal{P}_A$ is an open subset of the
af\/f\/ine space given by the intersection of $\operatorname{ker}(A)$ and the af\/f\/ine hull of $\Delta^{n-1}$, and
hence $\mathcal{P}_A$ is a~polytope of dimension $n-\ell-1$.

Under these hypotheses, we f\/irst establish Lemma~\ref{vertices}, demonstrating that the faces of $\mathcal{P}_A$
consist of the intersection of $\operatorname{ker}(A)$ with coordinate subspaces of $\mathbb{R}^n$.
If $I\subset\{1,\ldots,n\}$ is a~collection of indices, we again use the notation that
$V_I=\{(x_1,\ldots,x_n)\in\mathbb{R}^n\,|\, x_j=0
\ \forall\, j\notin I\}$ is the associated coordinate subspace.
We then show Proposition~\ref{simpcond}, which states Theorem~\ref{simpcondgen} for matrices that satisfy $(i)$, $(ii)$,
and $(iii)$, and then proceed to the proof of Theorem~\ref{simpcondgen}.
\begin{lemma}
\label{vertices}
Let $A\in\mathbb{Z}^{n\times\ell}$ satisfy hypotheses $(i)$, $(ii)$, and $(iii)$.
\begin{enumerate}[$(a)$]\itemsep=0pt
\item If $I\subset\{1,\ldots,n\}$ such that $\mathcal{P}_A\cap V_I=\{\nu\}$, then $\nu$ is a~vertex of
$\mathcal{P}_A$.
\item Each face $F$ of $\mathcal{P}_A$ is given by $F=\mathcal{P}_A\cap V_I$ for some $I$ with $|I|=\ell+\dim(F)+1$.
\end{enumerate}
\end{lemma}
As a~special case of $(b)$, note that each vertex of $\mathcal{P}_A$ is given by the intersection $\mathcal{P}_A\cap V_I$
where $I\subset\{1,\ldots,n\}$ is a~subset of cardinality $\ell+1$.
Note that given a~$k$-face $F$, the set $I$ given by condition $(b)$ need not be unique.
If $\operatorname{ker}(A)$ intersects the simplex $\Delta^{n-1}$ \emph{generically}, i.e., each of its vertices is
contained in the relative interior of an $\ell$-dimensional face of $\Delta^{n-1}$, then the $I$ corresponding to $F$
is unique.
In general, however, a~face can be contained in the intersection of several $(\ell+k)$-dimensional faces.
Given hypotheses $(i)$, $(ii)$, and $(iii)$, however, it is easy to see that a~vertex of $\mathcal{P}_A$ cannot correspond to
a~vertex of $\Delta^{n-1}$; this would indicate that a~standard basis vector $e_j$ is contained in the kernel, and
hence that the $j$th column of $A$ is a~zero column.
Similarly, if $I$ is a~set of indices of cardinality $|I|=\ell+1$, it need not be the case that $\mathcal{P}_A\cap V_I$
is a~singleton.
\begin{proof}
 $(a)$  Assume $\mathcal{P}_A\cap V_I=\{\nu\}$ for $I\subset\{1,\ldots,n\}$ with $\nu=(v_1,\ldots,v_n)$.
Suppose $\nu=tp+(1-t)q$ for $t\in]0,1[$ and $p=(p_1,\ldots,p_n),q=(q_1,\ldots,q_n)\in\mathcal{P}_A$.
For each $j\notin I$, we have that $tp_j+(1-t)q_j=v_j=0$.
As $p_j,q_j\geq0$, it follows that $p_j=q_j=0$.
Applying this argument to each $j\notin I$, it follows that $p,q\in V_I$.
Hence $p,q\in\mathcal{P}_A\cap V_I$, which was assumed to be a~singleton, so that $p=q=\nu$ and $\nu$ is a~vertex of~$\mathcal{P}_A$.

$(b)$ We prove the statement by induction on the codimension $c$ of the face $F$.
The case of $c=0$ is trivial.
Let $F$ be a~face of codimension $c+1$, so $k:=\dim(F)=n-\ell-c-2$.
Note that $F$ is contained in a~face $F'$ of codimension $c$.
By our inductive hypothesis, we can write $F'=\mathcal P_A\cap V_{I'}$ for some $I'\subset\{1,\dots,n\}$ of cardinality
$|I'|=\ell+(k+1)+1=n-c$.
This means that
\begin{gather*}
F'=\mathcal P_A\cap V_{I'}=\operatorname{ker}(A)\cap\Delta^{n+1}\cap V_{I'}=\operatorname{ker}(A)\cap\Delta^{n-c-1},
\end{gather*}
where $\Delta^{n-c-1}$ is the standard simplex in $V_{I'}$.

We claim that $F$ is contained in a~face of $\Delta^{n-c-1}$.
Letting $W:=\operatorname{af\/f}(\Delta^{n-c-1})\cap\operatorname{ker}(A)$ and $H_i^+:=\{(x_1,\dots,x_n)\in\mathbb
R^n|x_i\ge0\}$, we write
\begin{gather*}
F'=\operatorname{ker}(A)\cap\operatorname{af\/f}\big(\Delta^{n-c-1}\big)\cap\left(\cap_{i\in I'}H_i^+\right)=
W\cap\left(\cap_{i\in I'}H_i^+\right).
\end{gather*}
Setting $K_i^+:=W\cap H_i^+$, we have $F'=\cap_{i\in J}K_i^+$ for $J\subset I'$ chosen such that $K_i^+\ne W$ if and
only if $i\in J$.
It is a~well-known fact (see, e.g.,~\cite[Theorem 8.2]{BronstedBook}) that each facet of $F'$ is of the form $K_i\cap
W$ for some $i\in J$, where $K_i:=W\cap V_{\{i\}^c}$ is the supporting hyperplane of $K_i^+$.
Since~$F$ is a~facet of~$F'$, we conclude that
\begin{gather}
\label{Bronsted}
F=F'\cap K_i=F'\cap V_{\{i\}^c}\cap W.
\end{gather}
Since $F\subset F'\subset W$, it follows that that $F\subset F'\cap V_{\{i\}^c}$ which proves the claim.

Moreover, equation~\eqref{Bronsted} shows that $F=\mathcal P_A\cap V_I$ with $I:=I'\setminus\{i\}$.
\end{proof}
With this, we have the following.
\begin{proposition}
\label{simpcond}
Let $A$ be an $n\times\ell$ weight matrix satisfying hypotheses $(i)$, $(ii)$, and $(iii)$ so that~$\mathcal{P}_A$ is an
$(n-\ell-1)$-dimensional polytope.
Then conditions $(1)$, $(2)$, and $(3)$ of Theorem~{\rm \ref{simpcondgen}} are equivalent and are always
satisfied if $n\leq\ell+2$.
\end{proposition}
Note that given the hypotheses, condition (1) is equivalent to $\mathcal{P}_A$ having $n-\ell$ vertices.
Similarly, $A=A^\prime$ has full rank and no zero columns, simplifying~(2).
\begin{proof}
$(1) \Rightarrow (2)$: Suppose $\mathcal{P}_A$ has $n-\ell$ vertices $\nu_1,\ldots,\nu_{n-\ell}$.
To establish (2), we will show that each vertex $\nu_j$ lies in an $(\ell+1)$-dimensional coordinate plane, and the
intersection of these coordinate planes is an $\ell$-dimensional coordinate plane.
This will indicate the order of the vertices under which $A$ takes the required form.

For each vertex $\nu_j$, let $F_j$ denote the $(n-\ell-2)$-dimensional facet of $\mathcal{P}_A$ that does not contain~$\nu_j$, so that $F_j=\operatorname{cch}\{\nu_k\mid k\neq j\}$.
Then by Lemma~\ref{vertices}, each~$F_j$ is given by the intersection of~$\mathcal{P}_A$ with a~coordinate $n-1$-plane,
and hence corresponds to setting a~single coordinate equal to zero.
By reordering the variables $x_1,\ldots,x_n$, we may assume that $F_j=\mathcal{P}_A\cap V_{\{\ell+j\}^c}$ for
$j=1,\ldots,n-\ell$.
Let $v_{j,i}$ indicate the coordinates of $\nu_j$, i.e., $\nu_j=(v_{j,1},v_{j,2},\ldots,v_{j,n})$.
Note that for each $j$, as $V_{\{\ell+j\}^c}$ does not contain the vertex $\nu_j$, it follows that $v_{j,\ell+j}\neq0$.

For each $j$, we claim that $\cap_{k\neq j}F_k=\{\nu_j\}$.
To see this, f\/irst note that $\nu_j\in F_k$ for each $k\neq j$ so that $\{\nu_j\}\subset\cap_{k\neq j}F_k$.
For the reverse inclusion, suppose $p=(p_1,\ldots,p_n)\in\cap_{k\neq j}F_k$.
Then as $p\in\mathcal{P}_A$, we have that $p$ is a~convex combination of the $\nu_1,\ldots,\nu_{n-\ell}$, say
$p=\sum\limits_{m=1}^{n-\ell}t_m\nu_m$ with $0\leq t_m\leq1$ and $\sum\limits_{m=1}^{n-\ell}t_m=1$.
For each $r\neq j$, we have that $\cap_{k\neq j}F_k\subset F_r$ so that $p\in F_r$ and $p_{\ell+r}=0$.
As $v_{r,\ell+r}\neq0$, it then follows that $t_r=0$.
Therefore, the only nonzero $t_r$ is $t_j=1$, and $p=\nu_j$.
Letting $I_j=\{1,2,\ldots,\ell,\ell+j\}$, it follows that
\begin{gather*}
\{\nu_j\}=\bigcap\limits_{k\neq j}F_k=\bigcap\limits_{k\neq j}\big(\mathcal{P}_A\cap V_{\{\ell+k\}^c}\big)=
\mathcal{P}_A\cap\bigcap\limits_{k\neq j}V_{\{\ell+k\}^c}=\mathcal{P}_A\cap V_{I_j}.
\end{gather*}

Let $[D\mid C]$ denote the weight matrix $A$ row-reduced using integer scalar multiples, where $D$ is $\ell\times\ell$
and $C$ is $\ell\times(n-\ell)$.
We let $c_{k,j}$ denote the entries of $c$ as usual, with $1\leq k\leq\ell$ and $1\leq j\leq n-\ell$.
As $A$ has full rank, it must be that $[D\mid C]$ has full rank as well.
We claim that~$D$ is diagonal and nonsingular.

Suppose not, and then one of the pivot columns must be contained in $C$ so that the last row of $D$ is the zero row.
For each $j$, as $v_{j,\ell+k}=0$ for $k\neq j$, it follows that the $n$th entry of $[D\mid C]\nu_j$ is given by
$c_{\ell,j}v_{j,\ell+j}$.
Recall that $v_{j,\ell+j}\neq0$ and $\nu_j\in\operatorname{ker}(A)=\operatorname{ker}([D\mid C])$, and then
$c_{\ell,j}=0$.
However, as this is true for each $j\leq n-\ell$, it follows that the last row of $C$ is the zero row, contradicting
the fact that $[D\mid C]$ has full rank.
We conclude that $D$ is diagonal and nonsingular.
Clearly, by multiplying rows by~$-1$, we can assume that the diagonal entries of~$D$ are all negative.
Let~$d_k<0$ denote the diagonal entries of~$D$ for $1\leq k\leq\ell$.

Finally, we claim that each $c_{k,j}\geq0$.
For each $j$, as $\nu_j$ has nonzero coordinates only in the $1,2,\ldots,\ell$, and $\ell+j$ positions, we have that
the $k$th coordinate of $[D\mid C]\nu_j$ is given by $d_k v_{j,k}+c_{k,j}v_{j,\ell+j}$.
As $[D\mid C]\nu_j=0$, we have that $d_k v_{j,k}+c_{k,j}v_{j,\ell+j}=0$.
As $d_k<0$, as $v_{j,k}\geq0$, and as $v_{j,\ell+j}>0$.
It follows that $c_{k,j}\geq0$, completing the proof that~$(1)\Rightarrow(2)$.

$(2) \Rightarrow (3)$: Assuming $A$ is in the form $[D\mid C]$ as in (2), let~$\mu_i$ denote the
$i$th row of~$A$.
Then it is easy to see that $\langle\mu_i,e_i\rangle<0$ and $\langle\mu_i,e_j\rangle\geq0$ for $j\neq i$.
Moreover, $\operatorname{ker}(A)=\bigcap_{i=1}^\ell\mu_i^\perp$ by def\/inition.

\looseness=-1
$(3)\Rightarrow(1)$: Permute the coordinates $x_1,\ldots,x_n$ so that $j_i=i$ for
$i=1,\ldots,\ell$.
Let~$M$ be the $\ell\times n$ matrix with $i$th row $\mu_i$ and then $\operatorname{ker}(M)=\operatorname{ker}(A)$ by
hypothesis so that $\mathcal{P}_A=\mathcal{P}_M$.
Clearly, $M$~must then satisfy hypotheses $(i)$, $(ii)$, and $(iii)$, and moreover~$M$ is of the form $[D\mid C]$ as
described in condition (2).
Let $d_k<0$ denote the entries of $D$ and $c_{k,j}\geq0$ denote the entries of~$C$.

It is easy to see that each coordinate plane corresponding to $\{1,\ldots,\ell,\ell+k\}$ intersects $\mathcal{P}_M$ at
a~single vertex.
In particular, def\/ine
\begin{gather*}
\nu_j=\frac{1}{1+\sum\limits_{k=
1}^\ell-c_{k,j}/d_k}\left(\frac{-c_{1,j}}{d_1},\frac{-c_{2,j}}{d_2},\ldots,
\frac{-c_{\ell,j}}{d_\ell},0,\ldots,0,1,0,\ldots,0\right)
\end{gather*}
for $j=1,\ldots,n-\ell$, where the $1$ occurs in the $(\ell+j)$th position.
Simple computations show that each $\nu_j\in\operatorname{ker}([D\mid C])\cap\Delta^{n-1}$ and
$\operatorname{ker}([D\mid C])\cap V_{\{1,\ldots,\ell,\ell+j\}}$ is a~$1$-dimensional subspace of $\mathbb{R}^n$.
Therefore, $\{\nu_j\}=\operatorname{ker}([D\mid C])\cap\Delta^{n-1}\cap V_{\{1,\ldots,\ell,\ell+j\}}$, so that by
Lemma~\ref{vertices}, each $\nu_j$ is a~vertex of~$\mathcal{P}_M$.
It remains only to show that there are no other vertices.

However, for each $p=(p_1,\ldots,p_n)\in\mathcal{P}_M=\mathcal{P}_{[D\mid C]}$ the fact that $[D\mid C]p=0$ implies
that the $p_1,\ldots,p_\ell$ are uniquely determined by the $p_{\ell+1},\ldots,p_n$.
Moreover, letting $\pi\colon \mathbb{R}^n\to\mathbb{R}^{n-\ell}$ denote the projection
$\pi\colon (x_1,\ldots,x_n)\mapsto(x_{\ell+1},\ldots,x_n)$, it is obvious that
$\{\pi(\nu_1),\ldots,\pi(\nu_{n-\ell})\}$ is linearly independent in $\mathbb{R}^{n-\ell}$ and hence af\/f\/inely
independent.
Hence, given coordinates $p_{\ell+1},\ldots,p_n$, there is a~unique af\/f\/ine combination of the
$\pi(\nu_1),\ldots,\pi(\nu_{n-\ell})$ that yields $(p_{\ell+1},\ldots,p_n)$.
Then there are unique values $p_1,\ldots,p_\ell$ such that
$(p_1,\ldots,p_n)\in\operatorname{ker}(M)=\operatorname{ker}([D\mid C])$, and this af\/f\/ine combination of the
$\pi(\nu_j)$ is a~convex combination if and only $(p_1,\ldots,p_n)\in\Delta^{n-1}$.
It follows that each $p\in\mathcal{P}_M=\mathcal{P}_{[D\mid C]}$ is a~convex combination of the $\nu_j$, and hence that
there are no other vertices.
We conclude that the polytope $\mathcal{P}_A=\mathcal{P}_M=\mathcal{P}_{[D\mid C]}$ has $n-\ell$ vertices and hence, as
it is $(n-\ell-1)$-dimensional, that it is combinatorially equivalent to the standard $(n-\ell-1)$-simplex.

To complete the proof, we need only note that if $n\leq\ell+2$, then $\mathcal{P}_A$ is a~$0$- or $1$-dimensional
polytope, which is necessarily a~simplex.
\end{proof}

With this, we are prepared to prove Theorem~\ref{simpcondgen}, completing this subsection.
\begin{proof}
[Proof of Theorem~\ref{simpcondgen}] First, we note that the zero-f\/iber $J_A^{-1}(0)$ is contained in the preimage
under $\operatorname{sq}$ of the coordinate plane $V_{I_A}$ so that we may identify $\operatorname{sq}(J_A^{-1}(0))$
with $\operatorname{sq}(J_{A^\prime}^{-1}(0))$ via the embedding $\mathbb{R}^{|I_A|}\to\mathbb{R}^n$ induced by
$I_A\subset\{1,\ldots,n\}$.
Permute the coordinates $x_i$ for $i\in I_A$ so that any zero columns of $A^\prime$ are listed last.
Row reducing $A^\prime$ using integer scalar multiples yields a~matrix with any zero rows listed last of the form
\begin{gather*}
R=\left[
\begin{matrix}A^{\prime\prime}&0
\\
0&0
\end{matrix}
\right].
\end{gather*}
Here, $A^{\prime\prime}$ has dimensions $k\times(k+m)$ such that $k\leq\ell$ and $m\leq|I_A|-k$.
To see this, note that~$A^{\prime\prime}$ has a~pivot in each row by construction, and moreover that $A^{\prime\prime}$
has at least one positive and one negative element in each row to ensure that each~$x_i$ is nonzero for some element of
$\operatorname{ker}(A^{\prime\prime})$.

Clearly, the reduced space of the action of $\mathbb{T}^\ell$ on $\mathbb{R}^{|I_A|}$ with weight matrix $R$ coincides
with the reduced space of the action with weight matrix $R^\prime=[A^{\prime\prime}\, 0]$. Note that by construction,
$A^{\prime\prime}$~has full rank and no zero columns.
Moreover, for each $i\in I_A$, there is a~point in $\operatorname{ker}(A)$ with nonnegative coordinates such that
$x_i\geq0$.
It follows by convexity that $\operatorname{ker}(A)\cap(\Delta^{k+m-1})^\circ\neq\varnothing$.
Therefore, $A^{\prime\prime}$ satisf\/ies hypotheses $(i)$, $(ii)$, and $(iii)$.

Now,
\begin{gather*}
\mathcal{P}_{R^\prime}=\operatorname{ker}(R^\prime)\cap\Delta^{|I_A|-1}=
\big(\operatorname{ker}(A^{\prime\prime})\times\mathbb{R}^{|I_A|-m}\big)\cap\Delta^{|I_A|-1}
\\
\phantom{\mathcal{P}_{R^\prime}}{}
=\operatorname{cch}\big(\big(\operatorname{ker}\big(A^{\prime\prime}\big)\cap\Delta^{k+m-1}\big)\cup\{e_{m+1},\ldots,e_{|I_A|}\}\big),
\end{gather*}
where $e_1,\ldots,e_{|I_A|}$ denotes the standard basis of $\mathbb{R}^{|I_A|}$, $\Delta^{k+m-1}$ is the standard
simplex in $\mathbb{R}^{k+m}=\operatorname{Span}\{e_1,\ldots,e_{k+m}\}$, and the elements of
$\operatorname{ker}(A^{\prime\prime})\cap\Delta^{k+m-1}$ are identif\/ied with elements of $\mathbb{R}^{|I_A|}$ via the
obvious embedding $\mathbb{R}^{k+m}\to\mathbb{R}^{|I_A|}$.

With this, it is clear that the vertices of $\mathcal{P}_{R^\prime}$ are given by the $e_{m+1},\ldots,e_{|I_A|}$ along
with the images of the vertices of $\operatorname{ker}(A^{\prime\prime})\cap\Delta^{k+m-1}$ in $\mathbb{R}^{|I_A|}$ as
above.
Hence, $\mathcal{P}_{R^\prime}$ is a~polytope given by the closed convex hull of $\mathcal{P}_{A^{\prime\prime}}$ along
with $|I_A|-m$ points that are linearly independent to the vertices of $\mathcal{P}_{A^{\prime\prime}}$.
It follows that $A^\prime$ and hence $A$ is simplicial if and only if $A^{\prime\prime}$ is simplicial.
Recalling that $A^{\prime\prime}$ satisf\/ies hypotheses $(i)$, $(ii)$, and $(iii)$, an application of
Proposition~\ref{simpcond} to $A^{\prime\prime}$ completes the proof.
\end{proof}
\begin{example}
For the weight matrix given in Example~\ref{nonsimpexample}, the vertices of
$\mathcal{P}_A=\operatorname{ker}(A)\cap\Delta^5$ are given by $(1/2,1/2,0,0,0,0)$; $(1/2,0,0,1/2,0,0)$;
$(0,1/2,1/2,0,0,0)$; $(0,0,1/2,1/2,0,0)$; and $(0,0,0,0,1/2,1/2)$; so that $\mathcal{P}_A$ is a~$3$-dimensional
polytope with $5$ vertices.
Hence Theorem~\ref{simpcondgen}(1) fails, and $A$ is not simplicial.
\end{example}

\section{Smooth structures on singular phase spaces}
\label{smooth}

The aim of this section is to study singular phase spaces and smooth maps between them.
In this paper, we use the term \emph{singular phase space} loosely, i.e., we mean by it a~space (preferably with
singularities) on which one can do some sort of Hamiltonian mechanics.
In order to give precise def\/initions, there are some choices to be made.
It will be convenient for our purposes to focus on the notion of a~\emph{differential space} in the sense of
Sikorski~\cite{Sikorski}.

\subsection{Poisson dif\/ferential spaces}
\label{PoissonDF}
To begin, let us recall the def\/inition of a~stratif\/ied symplectic space and the theorem of Sjamaar and Lerman,
which says that every symplectic quotient is such a~space.
\begin{definition}
A \emph{stratified symplectic space} is a~Whitney stratif\/ied space $X=\sqcup_{i\in I}X_i$ with an algebra $\mathcal
C^\infty(X)$ of continuous functions such that
\begin{enumerate}[1)]\itemsep=0pt
\item each stratum $X_i$ is a~symplectic manifold, \item $\mathcal C^\infty(X)$ is a~Poisson algebra, and \item the
pullback $\mathcal C^\infty(X)\to\mathcal C^\infty(X_i)$ with respect to the inclusions $X_i\hookrightarrow X$ is
compatible with the Poisson bracket.
\end{enumerate}
The Poisson algebra $\mathcal C^\infty(X)$, is called the \emph{algebra of smooth functions} on $X$.
\end{definition}

If we regard $X$ merely as a~topological space, we say that $\mathcal C^\infty(X)$ is a~\emph{smooth structure} on~$X$.
In many cases (for example in the case of the theorem below), it is known (see e.g.~\cite{LermanSjamaar}) that one can reconstruct the stratif\/ication from the Poisson algebra $\mathcal C^\infty(X)$.
The question of when one can do so without using the Poisson bracket is, to our knowledge, open.
So morally, $\mathcal C^\infty(X)$ contains all the information about~$X$.
\begin{theorem}[\cite{LermanSjamaar}]\label{LSthm}
Let $G$ be a~compact Lie group acting on a~symplectic manifold $(M,\omega)$ in a~Hamiltonian way, and let
$J\colon  M\to\mathfrak g^*$ be a~moment map for this action.
Then the symplectic quotient $M_0=Z/G$, with $Z=J^{-1}(0)$, is a~stratified symplectic space, where the strata
\begin{gather*}
(M_0)_{(H)}:=(M_{(H)}\cap Z)/G
\end{gather*}
are indexed by conjugacy classes $(H)$ of subgroups $H\subset G$ that arise as isotropy groups of elements of~$Z$.
Here $M_{(H)}$ is the set of points in $M$ whose isotropy group is an element of the class~$(H)$.
The Poisson algebra of smooth functions $\mathcal C^\infty(M_0)$ is given by
\begin{gather*}
\mathcal C^\infty(M_0):=\mathcal C^\infty(M)^G/\big(\mathcal C^\infty(M)^G\cap I_Z\big),
\end{gather*}
where $I_Z\subset\mathcal C^\infty(M)$ denotes the ideal of smooth functions vanishing on $Z$, and $\mathcal
C^\infty(M)^G\subset C^\infty(M)$ is the Poisson subalgebra of of $G$-invariant smooth functions.
\end{theorem}

Note that elements of $\mathcal C^\infty(M_0)$ can be in fact regarded as functions on $M_0$.
Note further that it is not dif\/f\/icult to check that $\mathcal C^\infty(M)^G\cap I_Z\subset\mathcal C^\infty(M)$ is
actually a~Poisson ideal, so that the Poisson bracket on $\mathcal C^\infty(M_0)$ is canonically def\/ined.

Using the smooth structure as the key idea, one can easily talk about symplectomorphisms between stratif\/ied
symplectic spaces.

\begin{definition}[\cite{LermanSjamaarMontgomery}]\label{symplectomorphism}
A \emph{symplectomorphism} between the symplectic stratif\/ied spaces  $(X=\sqcup_{i\in I}X_i,\mathcal C^\infty(X))$ and $(Y=\sqcup_{j\in J}Y_j,\mathcal C^\infty(Y))$ is def\/ined to be a~homeomorphism
$\varphi\colon  X\to Y$ whose pullback $\mathcal C^\infty(Y)\to\mathcal C^\infty(X)$ with $f\mapsto
f\circ\varphi$ is an isomorphism of Poisson algebras.
\end{definition}

Before ref\/ining the concept of symplectomorphism, we recall the notion of a~Poisson dif\/fe\-ren\-tial
space~\cite{MultarzynskiZekanowski}.
This idea will help us strip of\/f the unnecessary details from the notion of a~stratif\/ied symplectic space and widen
the setup to including, e.g., orbit spaces of Poisson $G$-actions.
\begin{definition}
A \emph{differential space} (in the sense of Sikorski) is def\/ined as a~pair $(X,\mathcal C^\infty(X))$, where $X$ is
a~topological space and $\mathcal C^\infty(X)$ is an algebra of continuous functions on $X$ such that the following
axioms are fulf\/illed:
\begin{enumerate}\itemsep=0pt
\item The topology of $X$ is generated by $\mathcal C^\infty(X)$.
\item If $F\in\mathcal C^\infty(\mathbb R^n)$, $f_1,\dots,f_n\in\mathcal C^\infty(X)$, then
$F(f_1,\dots,f_n)\in\mathcal C^\infty(X)$.
\item If $f\colon  X\to\mathbb R$ is a~function such that for every $x\in X$, there exists an open
neighborhood $U$ of $x$ and an $f_U\in\mathcal C^\infty(X)$ such that $f_{|U}=f_U$, then $f\in\mathcal C^\infty(X)$.
\end{enumerate}
A \emph{Poisson differential space} is a~triple $(X,\mathcal C^\infty(X),\{\:,\:\})$, where $(X,\mathcal C^\infty(X))$
is a~dif\/ferentiable space and $\{\:,\:\}\colon \mathcal C^\infty(X)\times\mathcal C^\infty(X)\to\mathcal
C^\infty(X)$ is a~Poisson bracket.
\end{definition}
\begin{definition}
\label{def:smoothMap}
A \emph{smooth map} from the dif\/ferential space $(X,\mathcal C^\infty(X))$ to the dif\/ferential space $(Y,\mathcal
C^\infty(Y))$ is a~continuous map $\varphi\colon  X\to Y$ such that the pullback $\varphi^*\colon
f\mapsto f\circ\varphi$ sends smooth functions on $Y$ to smooth functions on $X$.
If in addition $\varphi^*\colon \mathcal C^\infty(Y)\to\mathcal C^\infty(X)$ preserves the Poisson structures,
$\varphi$ is called a~\emph{Poisson map}.
\end{definition}
If $(X,\mathcal C^\infty(X))$ is a~dif\/ferential space, then a~maximal ideal $\mathfrak m\subset\mathcal C^\infty(X)$
is called a~\emph{real maximal ideal} if the residue f\/ield $\mathcal C^\infty(X)/\mathfrak m$ is isomorphic to
$\mathbb R$.
The set $\operatorname{Spec}_{\mathbb R}(\mathcal C^\infty(X))$ of real maximal ideals in $\mathcal C^\infty(X)$ is
called the \emph{real spectrum} of $\mathcal C^\infty(X)$.

To make the def\/inition of a~Poisson dif\/ferential space $(X,\mathcal C^\infty(X),\{\:,\:\})$ workable, we have to
impose some additional assumptions, namely:
\begin{enumerate}[$(A)$]\itemsep=0pt
 \item The real spectrum consists of points, i.e., every real maximal ideal in $\mathcal C^\infty(X)$ is of the
form $\mathfrak m_\xi:=\{f\in\mathcal C^\infty(X)\mid f(\xi)=0\}$.
Moreover, we require $\bigcap_{\xi\in X}\mathfrak m_\xi=0$.
\item All Hamiltonian vector f\/ields $D$ (i.e., those of the form $D:=\{h,\:\}$ for some $h\in\mathcal C^\infty(X)$)
fulf\/ill the chain rule.
This means that, if we pick some $\varphi_1,\dots,\varphi_k\in\mathcal C^\infty(X)$ and put
$\varphi:=(\varphi_1,\dots,\varphi_k)\colon  X\to\mathbb R^k$, then for any $F\in\mathcal C^\infty(\mathbb
R^k)$ we have
\begin{gather*}
D(F\circ\varphi)=\sum_{i=1}^k\left(\frac{\partial F}{\partial x^i}\circ\varphi\right) D(\varphi_i).
\end{gather*}
\end{enumerate}

Note that condition $(A)$ is always satisf\/ied if $X$ is a~closed subset or $\mathbb{R}^n$ and $\mathcal{C}^\infty(X)$
is the quotient algebra $\mathcal{C}^\infty(\mathbb{R}^n)/I$ where $I$ is the closed ideal of
$\mathcal{C}^\infty(\mathbb{R}^n)$ consisting of functions that vanish on $X$; see~\cite[Proposition~2.13]{GonzalezBook}.
It is not clear to the authors if condition $(A)$ remains true, e.g.\
if~$X$ is not paracompact.
In addition, it is not known to the authors whether a~symplectic stratif\/ied space is automatically a~Poisson
dif\/ferential space fulf\/illing conditions~$(A)$ and~$(B)$.
However, using results from~\cite{CushmanSniatycki}, it is easy to prove the following.
\begin{proposition}
With the notation of Theorem~{\rm \ref{LSthm}}, if $M$ has a~finite number of orbit types as a~$G$-manifold, then the
symplectic quotient $(M_0,\mathcal C^\infty(M_0),\{\:,\:\})$ is a~Poisson differential space satisfying conditions
$(A)$ and $(B)$.
\end{proposition}
\begin{proof}
In~\cite{CushmanSniatycki} it is proven that for any action of a~compact Lie group $G$ on a~manifold~$M$, the space of
$G$-orbits $M/G$ is a~dif\/ferential space.
The smooth structure here is given by the the algebra $\mathcal C^\infty(M)^G$ of $G$-invariant functions on~$M$.
Every subspace of a~dif\/ferential space is a~dif\/ferential space, and hence the symplectic quotient is
a~dif\/ferential space.
Using a~system of dif\/ferentiable invariants (see Theorem~\ref{diffinv}), $M_0$ can be realized as a~closed
dif\/ferential subspace of $\mathbb{R}^k$ so that~\cite[Proposition 2.13]{GonzalezBook} applies.
Property $(A)$ follows.
Property $(B)$ is obvious.
\end{proof}

\subsection{Global charts and the lifting theorem}
\label{liftingthm}
In this subsection, we impose a~more rigid structure on our Poisson dif\/ferentiable spaces.
The terminology chosen stems from the observation that complete sets of dif\/ferentiable invariants (cf.\
the Schwarz--Mather Theorem~\ref{diffinv}) have much in common with linear coordinates on a~vector space.
In fact, those coordinates can be seen as a~Hilbert basis of a~trivial group representation.
\begin{definition}
\label{globalchart}
A \emph{global chart} on a~Poisson dif\/ferential space $(X,\mathcal C^\infty(X),\{\:,\:\})$ is an algebra homomorphism
\begin{gather*}
\varphi\colon \  \mathbb R[\boldsymbol x]:=\mathbb R[x_1,\dots,x_k]\to\mathcal C^\infty(X),
\qquad
x_i\mapsto\varphi_i,
\qquad
i\in\{1,\dots,k\},
\end{gather*}
such that
\begin{enumerate}\itemsep=0pt
\item The image of $\varphi$, denoted $\mathbb R[X]$, is a~Poisson subalgebra of $\mathcal C^\infty(X)$, called
\emph{Poisson subalgebra of regular functions} on $X$.
\item $\mathcal C^\infty(X)$ is $\mathcal C^\infty$-integral over $\mathbb R[\boldsymbol x]$, that is, for any
$f\in\mathcal C^\infty(X)$ there is a~$F\in\mathcal C^\infty(\mathbb R^k)$ such that $f=F\circ\varphi$.
Abusing language slightly, here $\varphi$ denotes the vector valued map $X\to\mathbb R^k$,
$\xi\mapsto(\varphi_1(\xi),\dots,\varphi_k(\xi))$.
\item The image of $\varphi$ in $\mathcal C^\infty(X)$ separates points.
\end{enumerate}

For a~global chart $\varphi\colon \mathbb R[\boldsymbol x]\to\mathcal C^\infty(X)$ we use property (1) to
transfer the Poisson structure from $\mathcal C^\infty(X)$ to $\mathbb R[\boldsymbol x]/\operatorname{ker}(\varphi)$.
In this way we obtain an embedding of Poisson algebras
\begin{gather*}%\label{PoissonEmb}
\overline{\varphi}\colon \ \mathbb R[\boldsymbol x]/\operatorname{ker}(\varphi)\hookrightarrow\mathcal C^\infty(X).
\end{gather*}
Of course, $\mathbb R[X]$ is isomorphic to $\mathbb R[\boldsymbol x]/\operatorname{ker}(\varphi)$.
If for the global chart $\varphi\colon \mathbb R[\boldsymbol x]\to\mathcal C^\infty(X)$, the algebra $\mathbb
R[\boldsymbol x]$ carries a~$\mathbb Z$-grading such that the ideal $\operatorname{ker}(\varphi)$ is homogeneous, we
call $\varphi\colon \mathbb R[\boldsymbol x]\to\mathcal C^\infty(X)$ a~\emph{$\mathbb Z$-graded global chart}.
\end{definition}

Our favorite examples are, of course, the symplectic quotients $M_0=Z/G$ (cf.\
Theorem~\ref{LSthm}).
In fact, if we pick a~complete system $\rho_1,\dots,\rho_k$ of dif\/ferentiable invariants for the $G$-action on~$M$,
see Theorem~\ref{diffinv}, then we can make out of it a~global chart by def\/ining $\varphi_i$ to be the class of
$\rho_i$ in $\mathcal C^\infty(M_0)=\mathcal C^\infty(M)^G/(\mathcal C^\infty(M)^G\cap I_Z)$ for each
$i\in\{1,\dots,k\}$ (confer Theorem~\ref{LSthm}).
Similarly, if we consider an orbit space of a~Poisson $G$-space with f\/initely many orbit types, we can take
a~complete system of invariants itself to form a~global chart.

Let us also comment on the linear case, which is our main concern in this paper.
If we examine the symplectic quotient coming from a~unitary representation $G\to\operatorname{U}_n$, then we should of
course choose a~minimal homogeneous system of polynomial invariants $\rho_1,\dots,\rho_k$.
If we assign to the variables $x_i$ in the above def\/inition the degree of $\rho_i$, the result is a~$\mathbb
Z$-graded global chart.
Similar comments apply to the orbit space of a~symplectic representation $G\to\operatorname{Sp}(\mathbb R^{2n})$.
Note that the choice of a~complete system of polynomial invariants, and therefore of a~global chart, is not unique.
This choice turns out not to be essential; see Remark~\ref{rem:HBasisChoice} below.
A~more severe problem is that it might be practically impossible to compute a~complete system of invariants.
As well, the determination of $\operatorname{ker}(\varphi)$ can be tricky.
\begin{lemma}
\label{phiisembedding}
With the notation of Definition~{\rm \ref{globalchart}}, the map $\varphi\colon  X\to\mathbb R^k$,
$\xi\mapsto(\varphi_1(\xi),\dots,\varphi_k(\xi))$ is injective.
\end{lemma}
\begin{proof}
By def\/inition, any regular function $f\in\mathbb R[X]\subset\mathcal C^\infty(X)$ can be written as the composition
$f=p\circ\varphi=p(\varphi_1,\dots,\varphi_k)$ of a~polynomial $p\in\mathbb R[\boldsymbol x]$.
So if $\varphi(\xi_1)=\varphi(\xi_2)$ for some $\xi_1,\xi_2\in X$, then $f(\xi_1)=f(\xi_2)$ for all $f\in\mathbb R[X]$.
By Def\/inition~\ref{globalchart}(3), $\mathbb R[X]$ separated points, and hence it follows that $\xi_1=\xi_2$.
\end{proof}

\begin{lemma}
\label{regulardense}
Assume that the Poisson differential space in Definition~{\rm \ref{globalchart}} has property $(A)$.
Then if $\epsilon\colon \mathcal C^\infty(X)\to\mathbb R$ is a~morphism of $\mathbb R$-algebras such that
$\epsilon_{|\mathbb R[X]}=0$ it follows that $\epsilon=0$.
\end{lemma}

\begin{proof}
Assume that $\epsilon$ is nonzero.
Then $\operatorname{ker}(\epsilon)$ is a~real maximal ideal and hence, by proper\-ty~$(A)$, of the form $\mathfrak
m_\xi=\{f\in\mathcal C^\infty(X)\mid f(\xi)=0\}$.
On the other hand, as $\mathbb R[X]$ separates points, there is an $f\in\mathbb R[X]$ such that $f(\xi)\ne0$.
This contradicts our assumption that $f\in\operatorname{ker}(\epsilon)$.
\end{proof}

We now def\/ine morphisms for Poisson dif\/ferential spaces with global charts.
\begin{definition}\label{arrow}
An \emph{arrow} from a~Poisson dif\/ferential space $(X,\mathcal C^\infty(X),\{\:,\:\})$ with global chart
$\varphi\colon \mathbb R[\boldsymbol x]=\mathbb R[x_1,\dots,x_k]\to\mathcal C^\infty(X)$ to a~Poisson
dif\/ferential space $(Y,\mathcal C^\infty(Y),\{\:,\:\})$ with global chart $\psi\colon \mathbb R[\boldsymbol
y]=\mathbb R[y_1,\dots,y_m]\to\mathcal C^\infty(Y)$ is a~morphism of algebras $\lambda\colon \mathbb
R[\boldsymbol y]\to\mathbb R[\boldsymbol x]$, such that
\begin{enumerate}[$(i)$]\itemsep=0pt \item We have $\lambda(\operatorname{ker}(\psi))\subset\operatorname{ker}(\varphi)$, and the induced morphism of
algebras
\begin{gather*}%\label{lambdabar}
\overline{\lambda}\colon\ \mathbb R[\boldsymbol y]/\operatorname{ker}(\psi)
\to\mathbb R[\boldsymbol x]/\operatorname{ker}(\varphi)
\end{gather*}
is compatible with the Poisson bracket.
\item Setting $\lambda_i:=\lambda(y_i)\in\mathbb R[\boldsymbol x]$, $i=1,\dots,m$, and def\/ining
\begin{gather*}
\vartheta\colon \  X\to\mathbb R^m,
\qquad
\vartheta(\xi):=\big((\varphi(\lambda_1))(\xi),\dots,(\varphi(\lambda_m))(\xi)\big),
\end{gather*}
the image $\operatorname{im}(\psi)$ of the map $\psi\colon  Y\to\mathbb R^m$ contains
$\operatorname{im}(\vartheta)$.
\end{enumerate}
If both charts are $\mathbb Z$-graded and the algebra morphism $\lambda$ is compatible with the grading we say that the
arrow is \emph{$\mathbb Z$-graded}.
\end{definition}

Clearly, an arrow contains redundant information~-- what is really important is $\overline{\lambda}$.
We say that two arrows $\lambda$ and $\lambda'$ are \emph{equivalent} if they induce the same $\overline{\lambda}$.
\begin{theorem}[lifting theorem]\label{lift}
With the notation of Definition~{\rm \ref{arrow}} and choice of an arrow $\lambda$, let us assume that both Poisson
differential spaces have property~$(B)$ and $(X,\mathcal C^\infty(X))$ has proper\-ty~$(A)$.
Then there exists unique morphism of Poisson algebras $\widetilde{\lambda}\colon \mathcal
C^\infty(Y)\to\mathcal C^\infty(X)$, such that $\varphi\circ\lambda=\widetilde{\lambda}\circ\psi$.
\begin{gather*}
\xymatrix{
\mathcal C^\infty(Y)\ar[r]^{\widetilde{\lambda}}&\mathcal C^\infty(X)\\
\mathbb R[\boldsymbol y]\ar[u]^\psi\ar[r]_\lambda&\mathbb R[\boldsymbol x]\ar[u]_\varphi}
\end{gather*}
Moreover, the \emph{lift} $\widetilde{\lambda}$ depends only on the equivalence class of $\lambda$ and can be
understood as the pullback of the continuous map
\begin{gather*}
\chi\colon \  X\to Y,
\qquad
\xi\mapsto\psi^{-1}(\vartheta(\xi)).
\end{gather*}
For two arrows $\lambda_1$ and $\lambda_2$ we have
$\widetilde{\lambda_1\circ\lambda_2}=\widetilde{\lambda_1}\circ\widetilde{\lambda_2}$.
\end{theorem}

Note that by construction, the map $\chi$ is smooth in the sense of Def\/inition~\ref{def:smoothMap}.
\begin{proof}
Take a~function $f\in\mathcal C^\infty(Y)$ and write it as a~composite with $\psi$,
\begin{gather*}
f(\eta)=F\big(\psi_1(\eta),\dots,\psi_m(\eta)\big)
\qquad
\forall\, \eta\in Y,
\end{gather*}
for some $F\in\mathcal C^\infty(\mathbb R^m)$.
The function $\widetilde{\lambda}(f)\in\mathcal C^\infty(X)$ is def\/ined to be
\begin{gather*}
\big(\widetilde{\lambda}(f)\big)(\xi):=F\big((\varphi(\lambda_i))(\xi),\dots,(\varphi(\lambda_m))(\xi)\big)=
F(\vartheta(\xi))
\qquad
\forall\, \xi\in X,
\end{gather*}
where $\lambda_i:=\lambda(y_i)\in\mathbb R[\boldsymbol x]$ for $i=1,\dots,m$.
Clearly $\widetilde{\lambda}$ does not depend the choice of $\lambda$ within its equivalence class and fulf\/ills
$\varphi\circ\lambda=\widetilde{\lambda}\circ\psi$.
By assumption (ii) of Def\/inition~\ref{arrow}, $\widetilde{\lambda}$ does not depend on the choice of $F$.

The uniqueness of $\widetilde{\lambda}$ is a~consequence of Lemma~\ref{regulardense}.
In fact, given another algebra morphism $\widehat{\lambda}\colon \mathcal C^\infty(Y)\to\mathcal C^\infty(X)$
such that $\varphi\circ\lambda=\widehat{\lambda}\circ\psi$, then
$\epsilon_\xi:=(\widetilde{\lambda}(f))(\xi)-(\widehat{\lambda}(f))(\xi)$ is an algebra morphism
$\epsilon_\xi\colon \mathcal C^\infty(X)\to\mathbb R$ whose restriction to $\mathbb R[X]$ vanishes.
So Lemma~\ref{regulardense} implies $\epsilon_\xi=0$ for all $\xi\in X$.
But this means that the function $\widetilde{\lambda}(f)-\widehat{\lambda}(f)$ vanishes everywhere on $X$, and is hence
zero by property (A).

By Def\/inition~\ref{arrow} and the injectivity of $\psi\colon  Y\to\mathbb R^k$ (cf.
Lemma~\ref{phiisembedding}), the map $\chi\colon \xi\mapsto\psi^{-1}(\vartheta(\xi))$ is well-def\/ined.
The verif\/ication of the claim $\chi^*=\widetilde{\lambda}$ is straightforward,
\begin{gather*}
\big(\chi^*f\big)(\xi)=f(\chi(\xi))=f\big(\psi^{-1}(\vartheta(\xi))\big)
=(F\circ\psi)\big(\psi^{-1}(\vartheta(\xi))\big)
=F(\vartheta(\xi))=\widetilde{\lambda}(f).
\end{gather*}

Finally, let us show that $\widetilde{\lambda}$ is compatible with the bracket.
By construction, for all $i,j\in\{1,\dots,m\}$ there is a~polynomial $\gamma_{ij}=\gamma_{ij}(y_1,\dots,y_m)\in\mathbb
R[\boldsymbol y]$ representing the class of $\{y_i,y_j\}$ in $\mathbb R[\boldsymbol y]/\operatorname{ker}(\psi)$.
We observe that
\begin{gather*}%\label{Yrel}
\{\psi_i,\psi_j\}=\{\psi(y_i),\psi(y_j)\}=\psi(\gamma_{ij}(y_1,\dots,y_m))
%\\
=\gamma_{ij}(\psi_1,\dots,\psi_m)\in\mathcal C^\infty(Y),
\end{gather*}
because $\psi$ is by def\/inition compatible with the bracket.
Since, by assumption, $\overline{\lambda}$ is compatible with the bracket, we see that
$\{\lambda_i,\lambda_j\}=\{\lambda(y_i),\lambda(y_j)\}\in\mathbb R[\boldsymbol x]$ coincides with
$\lambda(\gamma_{ij}(y_1,\dots,y_m))=\gamma_{ij}(\lambda_1,\dots,\lambda_m)\in\mathbb R[\boldsymbol x]$ up to
$\operatorname{ker}(\varphi)$.
It follows that
\begin{gather}
\label{Xrel}
\{\varphi(\lambda_i),\varphi(\lambda_j)\}=\varphi(\{\lambda_i,\lambda_j\})=
\gamma_{ij}(\varphi(\lambda_1),\dots,\varphi(\lambda_m))
=\gamma_{ij}\circ\vartheta\in\mathcal C^\infty(X).
\end{gather}
With these preparations, we compute for $\eta\in Y$, and $f=F\circ\psi$ and $g=G\circ\psi$, making use of property~$(B)$:
\begin{gather*}
\{f,g\}(\eta)
=\sum_{i,j=1}^m\frac{\partial F}{\partial x_i}\big(\psi_1(\eta),\dots,\psi_m(\eta)\big)
\frac{\partial G}{\partial x_j}\big(\psi_1(\eta),\dots,\psi_m(\eta)\big)\{\psi_i,\psi_j\}(\eta)
\\
\phantom{\{f,g\}(\eta)}
{}=\sum_{i,j=1}^m\frac{\partial F}{\partial x_i}\big(\psi_1(\eta),\ldots,\psi_m(\eta)\big)
\frac{\partial G}{\partial x_j}\big(\psi_1(\eta),\ldots,\psi_m(\eta)\big)
\gamma_{ij}\big(\psi_1(\eta),\dots,\psi_m(\eta)\big),
\end{gather*}
which yields for $\xi\in X$:
\begin{gather*}
\widetilde{\lambda}\left(\{f,g\}\right)(\xi)=
\sum_{i,j=1}^m\frac{\partial F}{\partial x_i}\big((\varphi(\lambda_1))(\xi),\dots,(\varphi(\lambda_m))(\xi)\big)
\\
\phantom{\widetilde{\lambda}\left(\{f,g\}\right)(\xi)=}{}
\times
\frac{\partial G}{\partial x_j}\big((\varphi(\lambda_1))(\xi),\dots,(\varphi(\lambda_m))(\xi)\big)
\gamma_{ij}\big((\varphi(\lambda_1))(\xi),\dots,(\varphi(\lambda_m))(\xi)\big)
\\
\phantom{\widetilde{\lambda}\left(\{f,g\}\right)(\xi)}
{}=\sum_{i,j=1}^m\frac{\partial F}{\partial x_i}(\vartheta(\xi))
\frac{\partial G}{\partial x_j}(\vartheta(\xi))\;\gamma_{ij}(\vartheta(\xi)).
\end{gather*}
On the other hand we have
\begin{gather*}
\{\widetilde{\lambda}(f),\widetilde{\lambda}(g)\}(\xi)
=\sum_{i,j=1}^m\frac{\partial F}{\partial x_i}\big((\varphi(\lambda_1))(\xi),\dots,(\varphi(\lambda_m))(\xi)\big)
\\
\phantom{\{\widetilde{\lambda}(f),\widetilde{\lambda}(g)\}(\xi)=}
{}\times
\frac{\partial G}{\partial x_j}\big((\varphi(\lambda_1))(\xi),\dots,(\varphi(\lambda_m))(\xi)\big)
 \{\varphi(\lambda_i),\varphi(\lambda_j)\}(\xi),
\end{gather*}
which in view of equation~\eqref{Xrel} implies that
$\widetilde{\lambda}\left(\{f,g\}\right)=\{\widetilde{\lambda}(f),\widetilde{\lambda}(g)\}$.
\end{proof}
\begin{definition}
A Poisson map $\chi$ between Poisson dif\/ferential spaces with global charts that is obtained as a~lift of an arrow
$\lambda$ as in Theorem~\ref{lift} is called a~\emph{regular Poisson map}.
If the arrow $\lambda$ is such that
\begin{enumerate}[1)]\itemsep=0pt
\item $\overline\lambda$ is an isomorphism, and \item (in the notation of Def\/inition~\ref{arrow})
$\operatorname{im}(\vartheta)=\operatorname{im}(\psi)$,
\end{enumerate}
then $\chi$ is called a~\emph{regular Poisson diffeomorphism}.
If the arrow is in addition $\mathbb Z$-graded we say that the regular Poisson map (resp.
regular Poisson dif\/feomorphism) is \emph{$\mathbb Z$-graded}.
\end{definition}

Regular Poisson dif\/feomorphisms between symplectic stratif\/ied spaces are examples of symplectomorphisms; see
Def\/inition~\ref{symplectomorphism}.
\begin{remark}
\label{rem:HBasisChoice}
Consider a~unitary representation $G\to\operatorname{U}_n$, let $M_0$ denote the associated symplectic quotient, and
let $\rho_1,\ldots,\rho_r$ and $\sigma_1,\ldots,\sigma_s$ denote two choices of minimal homogeneous systems of
polynomial invariants.
Let
\begin{gather*}
\varphi\colon \ \mathbb R[x_1,\dots,x_r]\to\mathcal C^\infty(M_0),
\qquad
x_i\mapsto\rho_i,
\qquad
i\in\{1,\dots,r\},
\end{gather*}
and
\begin{gather*}
\psi\colon \ \mathbb R[y_1,\dots,y_s]\to\mathcal C^\infty(M_0),
\qquad
y_i\mapsto\sigma_i,
\qquad
i\in\{1,\dots,s\},
\end{gather*}
denote the corresponding global charts for $M_0$.
Expressing each $\sigma_i$ in terms of the polynomials $\rho_1,\ldots,\rho_r$ def\/ines an arrow
$\lambda\colon \mathbb{R}[\boldsymbol{y}]\to\mathbb{R}[\boldsymbol{x}]$, and one checks that this is
a~$\mathbb{Z}$-graded regular Poisson dif\/feomorphism.
\end{remark}

In Section~\ref{counter}, we will argue that certain Poisson dif\/ferential spaces with global chart are not ($\mathbb
Z$-graded) regularly dif\/feomorphic, because their rings of regular functions are not isomorphic as ($\mathbb
Z$-graded) commutative $\mathbb R$-algebras.
This type of problem is entirely in the realm of the theory of commutative Noetherian rings.

There are many potential applications of the theory presented above, which we will indicate elsewhere.
Here, we will return to the consideration of toric symplectic quotients.

\subsection{Orbifold cases in dimension 2}
\label{magic}

The purpose of this subsection is twofold.
First of all, we will illustrate the machinery introduced in the last subsection by presenting a~concrete example of
a~$\mathbb Z$-graded regular symplectomorphism.
Secondly, we will show that the simplicial condition is actually suf\/f\/icient for a~two-dimensional symplectic
quotient to be symplectomorphic to an orbifold.
Note that by Theorem~\ref{simpcondgen}, any two-dimensional symplectic quotient corresponds to a~simplicial
representation.

Before doing so, we would like to comment on a~subtle point that one faces when determining the kernel of a~global
chart.
For example, if we are interested in the ideal of smooth functions on~$\mathbb C$ that vanish on the zero set of the
function $J=z\cc z$, it turns out that is generated not by~$J$ itself, but rather by the linear monomials~$z$ and~$\cc
z$.
In the next proposition, we indicate that we do not have to worry about this kind of problems in the situation at hand.
\begin{proposition}
\label{signchange}
Let $A\in\mathbb Z^{\ell\times n}$ be a~weight matrix that can, by elementary row operations and permutation of the
column indices, be brought into the form $A=[D\:|\:C]$, where $D\in\mathbb Z^{\ell\times\ell}$ is a~diagonal matrix
with strictly negative entries and $C\in\mathbb Z^{\ell\times(n-\ell)}$ has non-negative entries and no rows that are
identically zero.
Then the $G=\mathbb T^\ell$-invariant part $I^G_Z=I_Z\cap\mathcal\mathcal C^\infty(\mathbb C^n)^G$ of the vanishing
ideal $I_Z\subset\mathcal C^\infty(\mathbb C^n)$ is generated by the components $J_1,\dots,J_\ell\in C^\infty(\mathbb
C^n)^G$ of the moment map.
Here we view $I^G_Z$ as an ideal in $\mathcal C^\infty(\mathbb C^n)^G$.
\end{proposition}

\begin{proof}
Based on the signs of the entries of $A$, condition $(i)$ of~\cite[Proposition~2.2]{HerbigIyengarPflaum} is fulf\/illed,
and the result follows.
\end{proof}

Given the $G=\mathbb T^\ell$-action on $\mathbb C^n$ encoded by our weight matrix $A\in\mathbb Z^{\ell\times n}$, we
can f\/ind a~real Hilbert basis (i.e., complete set of real polynomial invariants) $\rho_1,\dots,\rho_k\in\mathbb
R[\mathbb C^n]^G$ such that $\rho_i=z_i\cc z_i$ for $i=1,\ldots,n$.
Because our group is abelian, the moment map itself is invariant.
We can express the components of the moment map in terms of the $\rho$'s using
\begin{gather}
\label{eq-Shell}
J_a=J_{e_a}=\frac{1}{2}\sum_{i=1}^n A_{ai}\rho_i,
\qquad
a=1,\dots,\ell.
\end{gather}
We will occasionally refer to the relations of the form $J_a=0$ as the \emph{shell relations}.
Furthermore, let us denote by $f_1,\dots,f_r\in\mathbb R[x_1,\dots,x_k]$ a~complete set of algebraic relations among
the $\rho_1,\dots\rho_k$.
Using this data, we construct a~global chart
\begin{gather*}
\varphi\colon \ \mathbb R[\boldsymbol x]=\mathbb R[x_1,\dots,x_k]\to\mathcal C^\infty(M_0),
\qquad
x_i\mapsto\varphi_i
\end{gather*}
for the symplectic quotient $M_0=J^{-1}(0)/\mathbb T^\ell$, where $\varphi_i$ is the image of $\rho_i$ in $\mathcal
C^\infty(M_0)$.
Proposition~\ref{signchange} enables us to determine the kernel of $\varphi$.
Then we have the following, which we expect remains true for an arbitrary weight matrix $A$ and will pursue this
elsewhere.
\begin{corollary}
\label{Kernel}
Under the assumptions of Proposition~{\rm \ref{signchange}}, the homogeneous ideal $\operatorname{ker}(\varphi)\subset\mathbb
R[\boldsymbol x]$ is the ideal generated by $f_1,\dots,f_r$ and the linear forms $g_a:=\sum\limits_{i=1}^n A_{ai} x_i$,
$a=1,\dots,\ell$.
\end{corollary}

\begin{proof}
From Proposition~\ref{signchange} it follows that the map $\mathbb
R[\mathbb{C}^n]^G\to\mathcal{C}^\infty(\mathbb{C}^n)^G$ gives rise to an injection
\begin{gather*}
\mathbb R[\mathbb{C}^n]^G/\langle J_1,\dots,J_\ell\rangle\to\mathcal{C}^\infty(\mathbb{C}^n)^G/I^G_Z.
\end{gather*}
Since $\varphi$ factors through this injection
\begin{gather*}
\xymatrix{\mathbb R[\boldsymbol x]\ar[rr]^{\varphi}\ar@{->>}[dr]&&\mathcal{C}^\infty(\mathbb{C}^n)^G/I^G_Z
\\
&\mathbb R[\mathbb{C}^n]^G/\langle J_1,\dots,J_\ell\rangle\ar[ur]&}
\end{gather*}
we have that $\operatorname{ker}(\varphi)$ is the kernel of the substitution homomorphism $\mathbb R[\boldsymbol
x]\to\mathbb R[\mathbb{C}^n]^G/\langle J_1,\dots,J_\ell\rangle$.
The claim now easily follows from the isomorphism theorems.
\end{proof}

We can assume without loss of generality that the weight matrix is of the form $A=[D\:|\:\boldsymbol{n}]$ where
$D=\operatorname{diag}(-a_1,\ldots,-a_\ell)$ is an $\ell\times\ell$ diagonal matrix with $a_1,\ldots,a_\ell>0$ and
$\boldsymbol{n}$ is a~single column with entries $n_1,\dots,n_\ell\geq0$.
We assume as well that the $G=\mathbb T^\ell$-action is ef\/fective, which implies that $\gcd(a_i,n_i)=1$ for each
$i\in\{1,2,\dots,\ell\}$.
Let us introduce the shorthand notation:
\begin{gather*}
\mathcal{A}:=\operatorname{lcm}(a_1,\dots,a_\ell),
\qquad
m_i:=\frac{n_i\mathcal{A}}{a_i}
\qquad
\mbox{for}
\quad
i=1,\dots,\ell,
\qquad
\mathcal{M}:=\sum_{i=1}^\ell m_i.
\end{gather*}
It is not dif\/f\/icult to show that
\begin{gather*}
\rho_{1}=\operatorname{Re}\left(z_{\ell+1}^{\mathcal{A}}\prod_{i=1}^\ell z_i^{m_i}\right),
\qquad
\rho_{2}=\operatorname{Im}\left(z_{\ell+1}^{\mathcal{A}}\prod_{i=1}^\ell z_i^{m_i}\right),
\qquad
\rho_{3}=z_{\ell+1}\cc z_{\ell+1},
\\
\rho_4=z_1\cc z_1,
\qquad
\dots,
\qquad
\rho_{\ell+3}=z_\ell\cc z_\ell.
\qquad
\end{gather*}
constitutes a~minimal real Hilbert basis of our $\mathbb T^\ell$-action on $\mathbb C^n$.
The degree of $\rho_1$ and $\rho_2$ is $\mathcal{A}+\mathcal{M}$, while the degree of $\rho_3,\dots,\rho_{\ell+3}$ is
two.
Using the language of the previous section, this leads to a~$\mathbb Z$-graded global chart
\begin{gather*}
\psi\colon \ \mathbb R[\boldsymbol y]=\mathbb R[y_1,\dots,y_{\ell+3}]\to\mathcal C^\infty(M_0),
\qquad
y_i\mapsto\psi_i
\end{gather*}
for our symplectic quotient $M_0=J^{-1}(0)/\mathbb T^\ell$, where $\psi_i$ is $\rho_i$ regarded as an element of
$\mathcal C^\infty(M_0)$.
The kernel $\operatorname{ker}(\psi)$ of the algebra morphism $\psi$ is generated by the polynomials
\begin{gather*}%\label{torusrelations}
y_1^2+y_2^2-\frac{\prod\limits_{i=1}^\ell m_i^{m_i}}{\mathcal{A}^\mathcal{M}} y_3^{\mathcal{A}+\mathcal{M}}
\qquad
\mbox{and}
\qquad
y_{3+i}-\frac{m_i}{\mathcal{A}} y_3
\qquad
\mbox{for}
\quad
i=1,\dots,\ell,
\end{gather*}
the latter coming from the shell relations (see equation~\eqref{eq-Shell}).
The image of the vector valued map
\begin{gather*}
\psi\colon \ M_0\to\mathbb R^{\ell+3},
\qquad
m\mapsto(\psi_1(m),\dots,\psi_{\ell+3}(m))
\end{gather*}
is determined by the semialgebraic condition
\begin{gather*}%\label{psiimage}
y_1^2+y_2^2-\frac{\prod\limits_{i=1}^\ell m_i^{m_i}}{\mathcal{A}^\mathcal{M}}  y_3^{\mathcal{A}+\mathcal{M}}=0,
\\
y_{3+i}-\frac{m_i}{\mathcal{A}} y_3=0
\qquad
\mbox{for} \quad i=1,\dots,\ell,
\qquad
y_3\ge0.
\end{gather*}

\begin{figure}[t]  \centering
\includegraphics[width=0.3\textwidth]{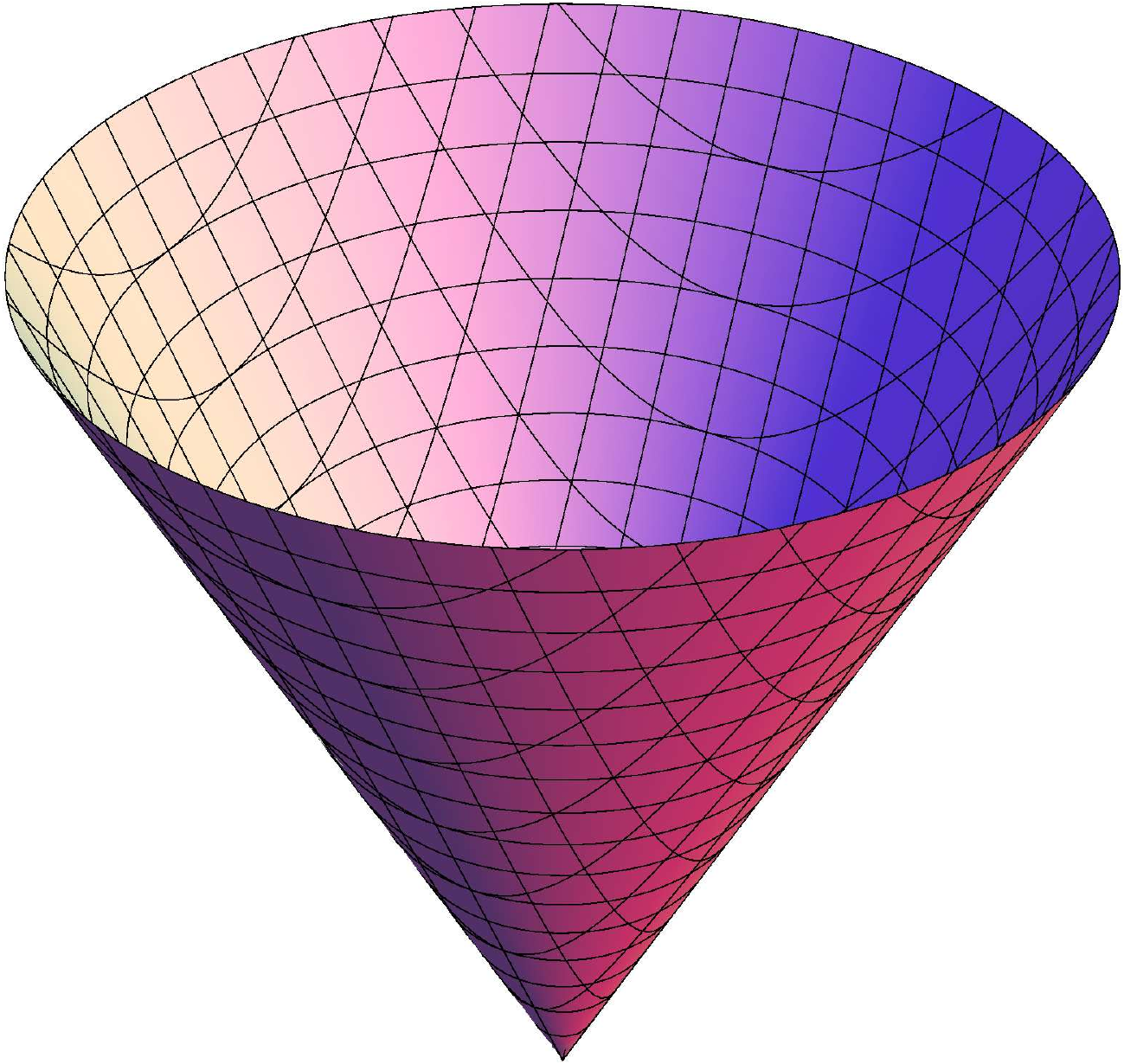} \qquad \quad
\includegraphics[width=0.3\textwidth]{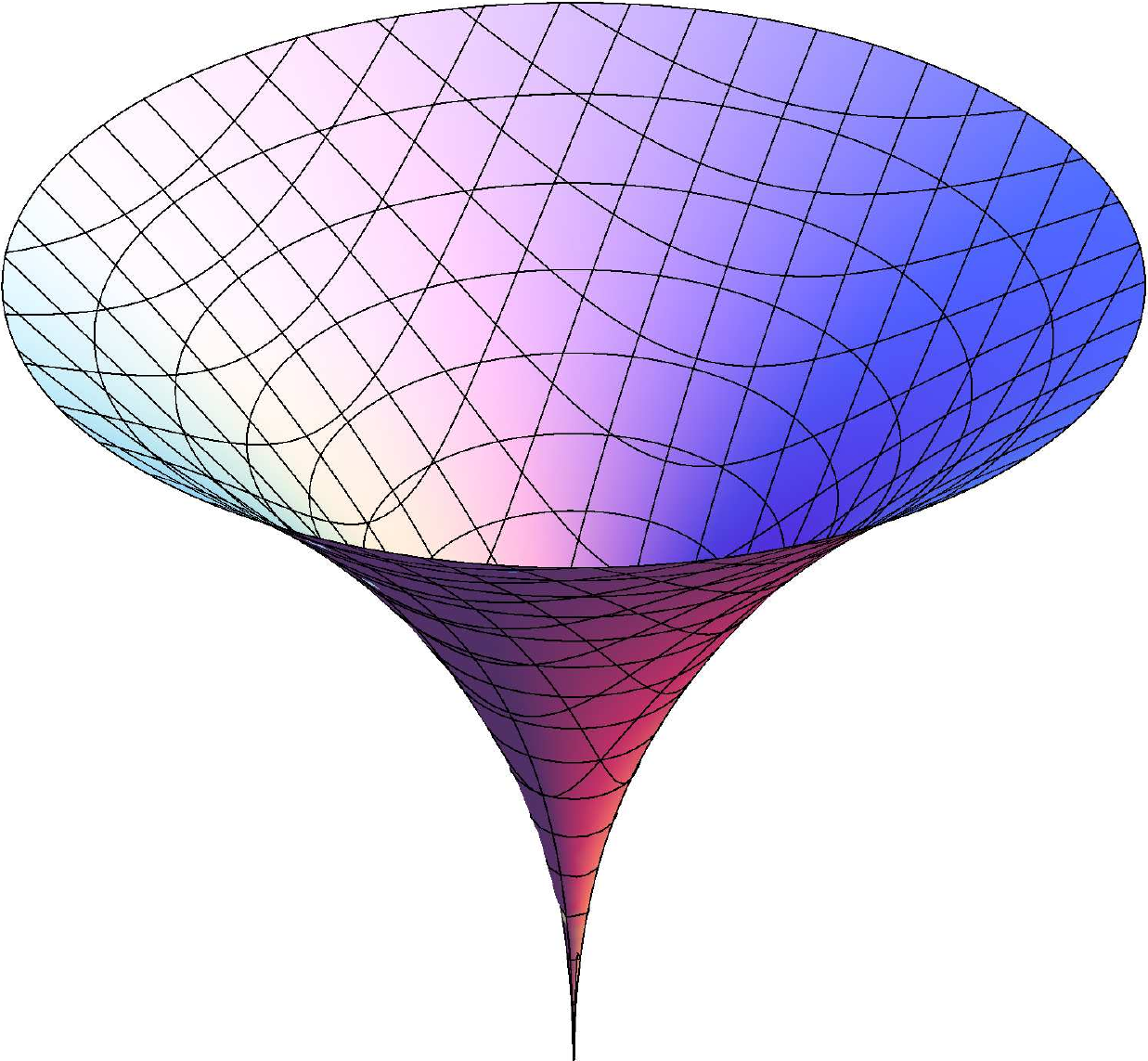}
\caption{The symplectic orbifold $\mathbb C/\mathbb Z_N$ for $N=2$ (left) and $N=5$ (right).}
\label{cones}
\end{figure}

On the other hand, let us consider the canonical action of the cyclic group $\mathbb Z_N$, for $N\ge2$, on $\mathbb C$.
In other words, we let $g\in\mathbb{Z}_N\subset\mathbb{S}^1\subset\mathbb C$ act on $z\in\mathbb C$ by multiplication.
Recall that the action on the complex conjugate variable $\cc z$ is given by $g^{-1}\cc z$.
As the $\mathbb Z_N$-action preserves the K{\"a}hler structure of $\mathbb C$, the quotient space $\left(X_N:=\mathbb
C/\mathbb Z_N,\mathcal C^\infty(X_N)=\mathcal C^\infty(\mathbb C)^{\mathbb Z_N}\right)$ is a~Poisson dif\/ferential
space of (real) dimension two.
It is easy to determine the real Hilbert basis consisting of $\varphi_1=\operatorname{Re}(z^N)$,
$\varphi_2=\operatorname{Im}(z^N)$, and $\varphi_3=z\cc z$.
Assigning to the variables $x_1$ and $x_2$ the degree $N$ and to $x_3$ the degree $2$, we obtain a~$\mathbb Z$-graded
global chart
\begin{gather*}
\varphi\colon \ \mathbb R[\boldsymbol x]=\mathbb R[x_1,x_2,x_3]\to\mathcal C^\infty(X_N),
\qquad
x_i\mapsto\varphi_i.
\end{gather*}
The kernel $\operatorname{ker}(\varphi)$ of the algebra morphism $\varphi$ is generated by the polynomial
$x_1^2+x_2^2-x_3^N$.
The image of the map $\varphi\colon  X_N\to\mathbb R^3$, $z\mapsto(\varphi_1(z),\varphi_2(z),\varphi_3(z))$,
is given by the semialgebraic set of solutions of the system
\begin{gather*}
x_1^2+x_2^2=x_3^N,
\qquad
x_3\ge0,
\end{gather*}
see Fig.~\ref{cones}.
With these preparations we are ready for the main result of this subsection.
\begin{theorem}
\label{magicthm}
With the above notation, if $N=\mathcal{A}+\mathcal{M}$, then the algebra homomorphism $\lambda\colon \mathbb
R[y_1,\dots,y_{\ell+3}]\longrightarrow\mathbb R[x_1,x_2,x_3]$ given by
\begin{gather*}
y_i\longmapsto\sqrt{\frac{\mathcal{A}^\mathcal{A}\prod\limits_{j=1}^\ell m_j^{m_j}}{N^N}} x_i,
\qquad
\text{for}
\quad
i=1,2,
\\
y_3\longmapsto\frac{\mathcal{A}}{N} x_3,
\qquad
y_{3+i}\longmapsto\frac{m_i}{N} x_3,
\qquad
\text{for}
\quad
i=1,\dots,\ell
\end{gather*}
is a~$\mathbb Z$-graded arrow lifting to a~$\mathbb Z$-graded symplectomorphism $X_N\to M_0$.
\end{theorem}
\begin{proof}
Using the relation $\{z_i,\cc z_j\}=-2\sqrt{-1}\delta_{ij}$, a~straightforward calculation yields
\begin{gather*}
\{\rho_1,\rho_2\}=
(z_{\ell+1}\cc z_{\ell+1})^{\mathcal{A}}\left(\prod_i(z_i\cc z_i)^{m_i}\right)
\left(\frac{\mathcal{A}^2}{z_{\ell+1}\cc z_{\ell+1}}+\sum_i\frac{m_i^2}{z_i\cc z_i}\right).
\end{gather*}
Writing $\cc y_i$ for the class of $y_i$ in $\mathbb R[\boldsymbol y]/\operatorname{ker}(\psi)$ and using $\cc
y_{i+3}=\frac{m_i}{\mathcal{A}}\cc y_3$ this leads to the relation
\begin{gather*}
\{\cc y_1,\cc y_2\}=
\frac{(\mathcal{A}+\mathcal{M})\prod_i m_i^{m_i}}{\mathcal{A}^{\mathcal{M}-1}}\cc y_3^{\mathcal{A}+\mathcal{M}-1}=
:\mathcal{B} \cc y_3^{\mathcal{A}+\mathcal{M}-1}.
\end{gather*}
Moreover, one can check that $\{\rho_1,\rho_3\}=2\mathcal{A}\rho_2$ and $\{\rho_2,\rho_3\}=-2\mathcal{A}\rho_1$.
We record our commutation relations $\{\cc y_i,\cc y_j\}$ in the table:
\begin{gather*}
\begin{array}{@{}c|ccc}&\cc y_1&\cc y_2&\cc y_3
\\
\hline
\cc y_1&0&\mathcal{B} \cc y_3^{\mathcal{A}+\mathcal{M}-1}&2\mathcal{A} \cc y_2
\phantom{\Big|}
\\
\cc y_2&&0&-2\mathcal{A} \cc y_1
\\
\cc y_3&&&0
\\
\end{array}
\end{gather*}
where we have omitted all $\cc y_{3+i}$, $i=1,\dots,\ell$ for the sake of brevity.

Similarly, we write $\cc x_i$ for the class of $x_i$ in $\mathbb R[\boldsymbol x]/\operatorname{ker}(\varphi)$.
We leave it to the reader to verify the multiplication table for the commutation relations $\{\cc x_i,\cc x_j\}$:
\begin{gather}
\begin{array}{@{}c|ccc}&\cc x_1&\cc x_2&\cc x_3
\\
\hline
\cc x_1&0&N^2 \cc x_3^{N-1}&2N \cc x_2
\phantom{\Big|}
\\
\cc x_2&&0&-2N \cc x_1
\\
\cc x_3&&&0
\\
\end{array}
\end{gather}
In order to construct the arrow $\lambda$ we make the def\/initions
\begin{gather*}
y_1\mapsto\alpha x_1,
\qquad
y_2\mapsto\alpha x_2,
\qquad
y_3\mapsto\beta x_3,
\end{gather*}
where $\alpha$, $\beta$ are determined from the multiplication tables, i.e.,
\begin{alignat*}{4}
& 2\alpha\beta N=2\mathcal{A} \alpha
\quad &&
\Rightarrow
\quad &&
\beta=\mathcal{A}/N, &
\\
& \alpha^2N^2=\beta^{N-1}\mathcal{B}
\quad &&
\Rightarrow
\quad &&
\alpha=\sqrt{\frac{(\mathcal{A}/N)^{N-1}\mathcal{B}}{N^2}}=\sqrt{\frac{\mathcal{A}^\mathcal{A}\prod\limits_{j=
1}^\ell m_j^{m_j}}{N^N}}. &
\end{alignat*}
Due to the identity $\alpha^2=\frac{\beta^N\mathcal{B}}{N\mathcal{A}}$, the generator
\begin{gather*}
y_1^2+y_2^2-\frac{\prod\limits_{i=1}^\ell m_i^{m_i}}{\mathcal{A}^\mathcal{M}} y_3^{\mathcal{A}+\mathcal{M}}=
y_1^2+y_2^2-\frac{\mathcal{B}}{N\mathcal{A}} y_3^{N}
\end{gather*}
of $\operatorname{ker}(\psi)$ is sent to the generator $x_1^2+x_2^2-x_3^N$ of $\operatorname{ker}(\varphi)$, which
proves that $\lambda$ is an arrow (the semialgebraic condition being clearly fulf\/illed).
By the lifting theorem, $\lambda$ lifts to a~Poisson map $X_N\to M_0$.
The inverse Poisson map can be constructed by lifting the arrow $x_1\mapsto\alpha^{-1}y_1$, $x_2\mapsto\alpha^{-1}y_2$
and $x_3\mapsto\beta^{-1}y_3$.
\end{proof}

Let us f\/inish this section by mentioning a~simple application of Theorem~\ref{magicthm}.
It is easy to see that if the weight matrix is of block form
\begin{gather*}
A=\left(
\begin{matrix}A_1&0
\\
0&A_2
\end{matrix}
\right),
\end{gather*}
then there is a~$\mathbb Z$-graded regular symplectomorphism from the symplectic quotient $M_A$ to \mbox{$M_{A_1}\times
M_{A_2}$}.
So if we consider for example (cf.\
\cite[p.~108]{HerbigIyengarPflaum}) the weight matrix $A\in\mathbb Z^{\ell\times2\ell}$ whose columns are given by $\pm e_i$,
where the $e_i$ are the standard basis vectors in $\mathbb R^\ell$, it follows that the reduced space $M_A$ is $\mathbb
Z$-graded regular symplectomorphic to the $\ell$-fold cartesian product of $\mathbb C/\mathbb Z_2$.

\section{Counterexamples in dimension~4}
\label{counter}

In this section, we prove that for certain unitary simplicial circle representations, there cannot exist a~$\mathbb
Z$-graded regular symplectomorphism from the symplectic quotient to a~quotient of a~linear symplectic action of
a~f\/inite group.
Before explaining our strategy, we note that a~natural idea is to examine ring theoretic features to distinguish
symplectic quotients from f\/inite quotients.
For example, it is known that invariant rings of unimodular representations of f\/inite groups are Gorenstein.
Unfortunately, cotangent lifted torus representations lead also to Gorenstein rings, because the representation
matrices are unimodular.
Since, provided the weight matrix $A$ has full rank, the shell relations cut out a~complete intersection in $\mathbb
R[\mathbb C^n]^G$, the rings of regular functions $\mathbb R[M_A]$ on the symplectic quotient space $M_A$ are
Gorenstein as well.

The only invariant we found useful in telling our symplectic quotients apart from f\/inite quotients as Poisson
dif\/ferential spaces is the \emph{Hilbert series} (also called the \emph{Poincar{\'e} series}) of the $\mathbb
Z$-graded ring of regular functions.
This is an invariant under $\mathbb{Z}$-graded regular symplectomorphism; whether it is an invariant under
symplectomorphism is not yet clear.
Let $V=\oplus_{i\ge0}V_i$ be a~positively graded, locally f\/inite-dimensional vector space over the f\/ield $\mathbb
K$.
Then the Hilbert series of $V$ is def\/ined as the formal power series
\begin{gather*}
\operatorname{Hilb}_{V|\mathbb K}(t)=\sum_{i\ge0}\dim_{\mathbb K}(V_i) t^i\in\mathbb Z[\![t]\!].
\end{gather*}
The Hilbert series of an invariant ring of a~compact groups can be calculated using Molien's formula (see e.g.~\cite{SturmfelsBook}).
It behaves well under cutting out complete intersections, which can be seen easily using the minimal free resolution.
Because our algebras of invariants are f\/initely generated, their Hilbert series can be written as $Q(t)/P(t)$, where
$Q(t)\in\mathbb Z[t]$ and $P(t)$ is of the form $\prod\limits_{i=1}^r(1-t^{n_i})^{k_i}$ with $k_i$ the number of generators in
degree $n_i$.
The Gorensteinness is ref\/lected by the fact that $Q(t)$ is palindromic.

Computations were performed using Singular\footnote{Decker~W., Greuel~G.M., Pf\/ister~G., Sch{\"o}nemann~H., {\sc Singular} {3-1-3}~-- {A}~computer algebra system for polynomial computations, 2011, \url{http://www.singular.uni-kl.de}.} %~\cite{DGPS}
and Mathematica\footnote{{Wolfram Research Inc.}, {\sc Mathematica} Edition: Version 7.0, Champaign, Illinois, 2008, \url{http://www.wolfram.com/mathematica/}}. %~\cite{Mathematica}.

\subsection[Weight matrices of type ${[}-1,1,m{]}$]{Weight matrices of type $\boldsymbol{[-1,1,m]}$}
%\label{11m}

Our f\/irst task here is to determine a~real Hilbert basis (i.e., complete systems of real homogeneous polynomial
invariants) for linear $G=\mathbb{S}^1$-actions on $\mathbb C^3$ corresponding to weight matrices of the form
$[-1,1,m]$ for $m=1,2,\dots$.
To this end we use the algorithm of Sturmfels~\cite{SturmfelsBook} together with the Groebner basis facilities of the
computer algebra system Singular.
There is one pitfall here, namely that we regard $\mathbb C^3$ as a~$6$-dimensional real vector space.
This means that $G$-operates also on the $\cc z$'s, so that we actually work with the weight matrix $[-1,1,m,1,-1,-m]$,
where the coordinates are ordered as $(z_1,z_2,z_3,\cc z_1,\cc z_2,\cc z_3)$.

The computer calculations indicate that the real Hilbert basis consists of the obvious real polynomials $z_1\cc z_1$,
$z_2\cc z_2$, and $z_3\cc z_3$, as well as the real and imaginary parts of
\begin{gather*}
z_1z_2
\qquad
\mbox{and}
\qquad
z_1^{m-i} \cc z^i_2z_3
\qquad
\mbox{for}
\quad
i=0,1,\dots,m,
\end{gather*}
giving altogether $7+2m$ polynomials.
This pattern has been verif\/ied for $m=1,\ldots,10$ using Singular.
The next task is to determine the algebraic relations among these generators.
We restrict to the cases $m=1$ and $2$; the rest of this subsection is devoted to a~more detailed presentation of these
two cases.

$\boldsymbol{m=1.}$
Here, the Hilbert basis consists of $9$ polynomials of degree $2$:
\begin{gather*}
\rho_1=z_1\cc z_1,
\qquad
\rho_2=z_2\cc z_2,
\qquad
\rho_3=z_3\cc z_3,
\qquad
\rho_4=\operatorname{Re}(z_1z_2),
\qquad
\rho_5=\operatorname{Im}(z_1z_2),
\\
\rho_6=\operatorname{Re}(z_1z_3),
\qquad
\rho_7=\operatorname{Im}(z_1z_3),
\qquad
\rho_8=\operatorname{Re}(z_2\cc z_3),
\qquad
\rho_9=\operatorname{Im}(z_2\cc z_3).
\end{gather*}
Among them we have $9$ quadratic relations (seen, of course, as being of degree 4):
\begin{alignat*}{4}
& \rho_4^2+\rho_5^2-\rho_1\rho_2=0,\qquad && \rho_3\rho_4+\rho_7\rho_9-\rho_6\rho_8=0,\qquad && \rho_6\rho_9+\rho_7\rho_8-\rho_3\rho_5=0,&
\\
& \rho_6^2+\rho_7^2-\rho_1\rho_3=0,\qquad && \rho_4\rho_8+\rho_5\rho_9-\rho_2\rho_6=0,\qquad && \rho_2\rho_7+\rho_4\rho_9-\rho_5\rho_8=0,&
\\
& \rho_8^2+\rho_9^2-\rho_2\rho_3=0,\qquad && \rho_4\rho_6+\rho_5\rho_7-\rho_1\rho_8=0,\qquad && \rho_4\rho_7+\rho_5\rho_6-\rho_1\rho_9=0.&
\end{alignat*}
That the numbers of generators and relations coincide here seems to be an accident.
Using the language of Subsection~\ref{lift}, we construct a~$\mathbb Z$-graded global chart $\varphi:\mathbb
R[\boldsymbol x]=\mathbb R[x_1,x_2,\dots,x_9]\to\mathcal C^{\infty}(M_0)$, of the form $x_i\mapsto\varphi_i$, where
$\varphi_i$ is simply $\rho_i$ seen as an element of $\mathcal C^{\infty}(M_0)$, and we assign to all~$x_i$ the degree~$2$.
In order to determine $\operatorname{ker}(\varphi)$, we use the fact that the representation is simplicial and
therefore, by Proposition~\ref{Kernel}, the only additional relation comes from the shell relation
$J=(\rho_2+\rho_3-\rho_1)/2=0$.
Summing up, we conclude that the $\mathbb Z$-graded ring $\mathbb R[M_0]$ of regular functions on the symplectic
quotient $M_0=Z/G$ is isomorphic to the ring $\mathbb R[x_1,x_2,\dots,x_9]/I$, where $I=\operatorname{ker}(\varphi)$ is
the homogeneous ideal
\begin{gather*}
I=\langle x_4^2+x_5^2-x_1x_2,\  x_6^2+x_7^2-x_1x_3,\  x_8^2+x_9^2-x_2x_3,\  x_3x_4+x_7x_9-x_6x_8,
\\
\hphantom{I=\langle}
x_4x_8+x_5x_9-x_2x_6,\  x_2x_7+x_4x_9-x_5x_8,\  x_6x_9+x_7x_8-x_3x_5,
\\
\hphantom{I=\langle}
x_4x_6+x_5x_7-x_1x_8,\  x_4x_7+x_5x_6-x_1x_9,\  x_2+x_3-x_1\rangle.
\end{gather*}
It is convenient to record some of the information contained in the minimal free resolution of the $\mathbb Z$-graded
$\mathbb R[x_1,x_2,\dots,x_9]$-module $\mathbb R[x_1,x_2,\dots,x_9]/I$ in the so-called \emph{Betti table} (for more
details see, e.g.,~\cite{Eisenbud}).
For the example at hand, we have computed the Betti table using Singular; see Table~\ref{BTab1}.
\begin{table}[t] \scriptsize
\centering
\caption{The Betti table for the symplectic quotient associated to $[-1,1,1]$.}\label{BTab1}
\vspace{1mm}

\begin{tabular}{c|cccccc}
& 0 & 1 & 2 & 3 & 4 & 5
\\
\hline
0 & 1 & -- & -- & -- & -- & --
\\
1 & -- & 1 & -- & -- & -- & --
\\
2 & -- & -- & -- & -- & -- & --
\\
3 & -- & 9 & -- & -- & -- & --
\\
4 & -- & -- & 25 & -- & -- & --
\\
5 & -- & -- & -- & 25 & -- & --
\\
6 & -- & -- & -- & -- & 9 & --
\\
7 & -- & -- & -- & -- & -- & --
\\
8 & -- & -- & -- & -- & 1 & --
\\
9 & -- & -- & -- & -- & -- & 1
\\
\hline
\mbox{total}& 1 & 10 & 25 & 25 & 10 & 1
\end{tabular}
\end{table}

It is easy to read of\/f from the Betti table the Hilbert series of $\mathbb R[M_0]$:
\begin{gather}
\nonumber
\operatorname{Hilb}_{\mathbb R[M_0]|\mathbb R}(t)=\frac{1-t^2-9t^4+25t^6-25t^8+9t^{10}+t^{12}-t^{14}}{(1-t^2)^9}
=\frac{1+4t^2+t^4}{(1-t^2)^4}
\\
\phantom{\operatorname{Hilb}_{\mathbb R[M_0]|\mathbb R}(t)}
{}=1+8t^2+27t^4+64t^6+125t^8+\cdots
\stackrel{*}{=}\sum\limits_{n=0}^\infty(n+1)^3t^{2n}.
\label{Hilb111}
\end{gather}
The step $(*)$ follows easily from the identity $\sum_i(-1)^i{4\choose i}(k+i)^3=0$, $k\ge0$,
which in turn can be proved by induction.
\hfill $\triangle$

$\boldsymbol{m=2.}$ Here the real Hilbert basis consists of eleven elements:
\begin{gather*}
\mbox{degree2:} \left\{
\begin{array}{@{}lll}
\rho_1=z_1\cc z_1,&\rho_2=z_2\cc z_2,&\rho_3=z_3\cc z_3,
\\
\rho_4=\operatorname{Re}(z_1z_2),&\rho_5=\operatorname{Im}(z_1z_2),
\end{array}\right.
\\
\mbox{degree3:}
\left\{
\begin{array}{@{}lll}
\rho_6=\operatorname{Re}\big(z_1^2z_3\big),
&
\rho_8=\operatorname{Re}(z_1\cc z_2z_3),
&
\rho_{10}=\operatorname{Re}(z_2^2\cc z_3),
\\
\rho_7=\operatorname{Im}\big(z_1^2z_3\big),
&
\rho_9=\operatorname{Im}(z_1\cc z_2z_3),
&
\rho_{11}=\operatorname{Im}(z_2^2\cc z_3).
\end{array}\right.
\end{gather*}
According to Singular, we have altogether 24 relations:
\begin{gather*}
\text{degree4:}\left\{\rho_4^2+\rho_5^2-\rho_1\rho_2=0,
\right.
\\
\text{degree5:}\left\{
\begin{array}{@{}ll}\rho_4\rho_{10}+\rho_5\rho_{11}-\rho_2\rho_8=0,&\rho_2\rho_9+\rho_4\rho_{11}-\rho_5\rho_{10}=0,
\\
\rho_5\rho_{9}+\rho_2\rho_{6}-\rho_4\rho_8=0,&\rho_4\rho_{9}+\rho_5\rho_{8}-\rho_2\rho_7=0,
\\
\rho_4\rho_{8}+\rho_5\rho_{9}-\rho_1\rho_{10}=0,&\rho_4\rho_{9}+\rho_1\rho_{11}-\rho_5\rho_8=0,
\\
\rho_4\rho_{6}+\rho_5\rho_{7}-\rho_1\rho_{8}=0,&\rho_1\rho_{9}+\rho_5\rho_{6}-\rho_4\rho_7=0,
\end{array}
\right.
\\
\text{degree6:}\left\{
\begin{array}{@{}ll}\rho_{9}\rho_{11}-\rho_{8}\rho_{10}+\rho_{2}\rho_{3}\rho_{4}=
0,&\rho_{8}\rho_{11}+\rho_{9}\rho_{10}-\rho_{2}\rho_{3}\rho_{5}=0,
\\
\rho_{6}\rho_{8}+\rho_{7}\rho_{9}-\rho_{1}\rho_{3}\rho_{4}=
0,&\rho_{6}\rho_{9}-\rho_{7}\rho_{8}+\rho_{1}\rho_{3}\rho_{5}=0,
\\
\rho_{6}^2+\rho_{7}^2-\rho_{1}^2\rho_{3}=0,&\rho_{3}(\rho_{4}^2-\rho_{5}^2)+\rho_{7}\rho_{11}-\rho_{6}\rho_{10}=0,
\\
\rho_{8}^2+\rho_{9}^2-\rho_{1}\rho_{2}\rho_{3}=0,&\rho_{9}^2-\rho_{8}^2+\rho_{6}\rho_{10}+\rho_{7}\rho_{11}=0,
\\
2\rho_{8}\rho_{9}+\rho_{6}\rho_{11}-\rho_{7}\rho_{10}=0,&\rho_{6}\rho_{11}+\rho_{7}\rho_{10}-2\rho_{3}\rho_{4}\rho_{5}=
0,
\\
\rho_{10}^2+\rho_{11}^2-\rho_2^2\rho_3=0,
\end{array}
\right.
\\
\text{degree8:}\left\{
\begin{array}{@{}l}\rho_1\big(\rho_{10}^2+\rho_{11}^2\big)-\rho_2\big(\rho_{8}^2+\rho_{9}^2\big)=0,
\\
\rho_1\big(\rho_{8}^2+\rho_{9}^2\big)-\rho_2\big(\rho_{6}^2+\rho_{7}^2\big)=0,
\\
\rho_1(\rho_{8}\rho_{10}-\rho_{9}\rho_{11})-\rho_2(\rho_{6}\rho_{8}+\rho_{7}\rho_{9})=0,
\\
\rho_1(\rho_{8}\rho_{11}+\rho_{9}\rho_{10})+\rho_2(\rho_{6}\rho_{9}-\rho_{7}\rho_{8})=0.
\end{array}
\right.
\end{gather*}
We would like to point out that it is, in principle, often more convenient to work with the complexif\/ication of our
base ring $\mathbb R[\mathbb C^n]$.
Doing so, we can f\/ind a~Hilbert basis consisting of monomials, e.g., $z_1z_2=\rho_4+\sqrt{-1}\rho_5$, $\cc z_1\cc
z_2=\rho_4-\sqrt{-1}\rho_5$ etc.
As a~consequence of complexif\/ication, the relations can be written as binomials.

Analogous to the case $m=1$, we obtain a~$\mathbb Z$-graded global chart for the symplectic quotient $M_0=Z/G$,
\begin{gather*}
\varphi: \ \mathbb R[\boldsymbol x]=\mathbb R[x_1,x_2,\dots,x_{11}]\to\mathcal C^{\infty}(M_0),
\qquad
x_i\mapsto\varphi_i,
\end{gather*}
where $\varphi_i$ is $\rho_i$ seen as an element in $\mathcal C^{\infty}(M_0)$, and we assign to $x_1,\dots,x_5$ the
degree $2$ and to $x_6,\dots,x_{11}$ the degree $3$.
Taking into account the shell relation $J=(\rho_2+2\rho_3-\rho_1)/2=0$, we see that the homogeneous ideal
$I=\operatorname{ker}(\varphi)\subset\mathbb R[\boldsymbol x]$ is generated by $25$ polynomials.
For the sake of brevity we leave it to the reader to write them down.
Once again, Singular is able to compute the minimal free resolution of the $\mathbb R[\boldsymbol x]$-module $\mathbb
R[\boldsymbol x]/I$.
The Betti table given in Table~\ref{BTab2}.
\begin{table}[t] \scriptsize \centering
\caption{The Betti table for the symplectic quotient associated to $[-1,1,2]$.}\label{BTab2}
\vspace{1mm}

\begin{tabular}{c|cccccccc} & 0& 1& 2& 3& 4& 5& 6& 7
\\
\hline
0& 1& --& --& --& --& --& --& --
\\
1& --& 1& --& --& --& --& --& --
\\
2& --& --& --& --& --& --& --& --
\\
3& --& 1& --& --& --& --& --& --
\\
4& --& 8& 1& --& --& --& --& --
\\
5& --& 11& 16& --& --& --& --& --
\\
6& --& --& 43& 8& --& --& --& --
\\
7& --& --& 24& 53& --& --& --& --
\\
8& --& --& --& 72& 21& --& --& --
\\
9& --& --& --& 21& 72& --& --& --
\\
10& --& --& --& --& 53& 24& --& --
\\
11& --& --& --& --& 8& 43& --& --
\\
12& --& --& --& --& --& 16& 11& --
\\
13& --& --& --& --& --& 1& 8& --
\\
14& --& --& --& --& --& --& 1& --
\\
15& --& --& --& --& --& --& --& --
\\
16& --& --& --& --& --& --& 1& --
\\
17& --& --& --& --& --& --& --& 1
\\
\hline
\text{total}& 1& 21& 84& 154& 154& 84& 21& 1
\end{tabular}
\end{table}
From this it is easy to derive a~formula for the Hilbert series $\operatorname{Hilb}_{\mathbb R[M_0]|\mathbb
R}(t)=Q(t)/P(t)$:
\begin{gather*}
Q(t)=1-t^2-t^4-8t^5-10t^6+16t^7+43t^8+16t^9-53t^{10}-72t^{11}+72t^{13}
\\
\phantom{Q(t)=}
{} +53t^{14}-16t^{15}-43t^{16}-16t^{17}+10t^{18}+8t^{19}+t^{20}+t^{22}-t^{24},
\\
P(t)=\big(1-t^2\big)^5\big(1-t^3\big)^6.
\end{gather*}
This can be reduced to
\begin{gather*}
\frac{1+2t^2+4t^3+2t^4+t^6}{(1-t^2)^2(1-t^3)^2}=1+4t^2+6t^3+9t^4+16t^5+26t^6+30t^7+\cdots.
\tag*{$\triangle$}
\end{gather*}

Let us close this subsection with the simple observation that the GIT-quotient $\mathbb C^3/\!/\mathbb{T}_\mathbb{C}^1$
corresponding to weight matrix $A=[-1,1,m]$ is the af\/f\/ine space $\mathbb C^2$.
In fact, $t\in\mathbb{T}_\mathbb{C}^1$ acts on $(q_1,q_2,q_3)\in\mathbb C^3$ by $t.q_1=t^{-1}q_1$, $t.q_2=tq_2$ and
$t.q_3=t^mq_3$.
A complex Hilbert basis is provided by the algebraically independent polynomials $q_1q_2,q_1^m q_3\in\mathbb
C[q_1,q_2,q_3]^{\mathbb{T}_\mathbb{C}^1}=\mathbb C[\mathbb C^3]^{\mathbb{T}_\mathbb{C}^1}$.

\subsection[Finite subgroups of ${\rm U}_2$]{Finite subgroups of $\boldsymbol{\operatorname{U}_2}$}
\label{Coxeter}

In this subsection, we consider the f\/inite subgroups $G$ of $\operatorname{U}_2$.
The classif\/ication of such subgroups is due to P.~du Val~\cite{Coxeter, DuVal}, and the details of the classif\/ication are recalled in Appendix~\ref{CoxApp}.
Of particular importance is the ADE-classif\/ication of f\/inite subgroups of $\operatorname{SU}_2$.
Here, the quotients $\mathbb C^2/G$ lead to the famous Kleinian singularities, while quotients $(\mathbb
C^2\times\overline{\mathbb C^2})/G$ appear not to be as well-understood.

Rather than an exhaustive computation of the Hilbert series for the nine families of f\/inite subgroups
$\operatorname{U}_2$, it will be suf\/f\/icient for our purposes to determine the Hilbert series for certain subgroups
$G<\operatorname{SU}_2$.
By the ADE-classif\/ication, any such group is cyclic or conjugate to a~binary dihedral, tetrahedral, octahedral, or
icosahedral group.

Recall the formula of Molien~\cite{Molien}, see e.g.~\cite[Theorem~2.2.1]{SturmfelsBook}, which expresses the Hilbert series for the $G$-invariant polynomials on
$\mathbb{C}^2$ as
\begin{gather}
\label{MolSeries}
\operatorname{Hilb}_{\mathbb{R}[\mathbb{C}^2]^G|\mathbb{R}}(t)=
\frac{1}{|G|}\sum\limits_{g\in G}\frac{1}{\det(\operatorname{id}-tg)}.
\end{gather}
In order to evaluate these summations, we will make use of a~Dedekind sum formula of I.~Gessel~\cite{Gessel}.

\subsubsection{Cyclic groups}
\label{subsubsec:CyclicGps}

Suppose $G\cong\mathbb{Z}_N$ is a~cyclic group of order $N$ with generator
$\operatorname{diag}(\omega_N,\omega_N^{-1})$ in complex coordinates, where $\omega_N$ is a~primitive $N$th root of
unity.
The action of this generator on $\mathbb C^2\times\overline{\mathbb C^2}$ is given by the matrix
$\alpha_N=\operatorname{diag}(\omega_N,\omega_N^{-1},\omega_N^{-1},\omega_N)$.
We will give explicit computations of the Hilbert series for two cases of particular interest separately followed by
the general case.

$\boldsymbol{N=2.}$ In this case, as $\omega_2=\omega_2^{-1}$, all quadratic polynomials are
invariant, and a~Hilbert basis is given by
\begin{gather*}
\rho_1=z_1\cc z_1,
\qquad
\rho_2=z_2\cc z_2,
\qquad
\rho_3=\operatorname{Re}(z_1z_2),
\qquad
\rho_4=\operatorname{Im}(z_1z_2),
\qquad
\rho_5=\operatorname{Re}(z_1\cc z_2),
\\
\rho_6=\operatorname{Im}(z_1\cc z_2),
\qquad
\rho_7=\operatorname{Re}\big(z_1^2\big),
\qquad
\rho_8=\operatorname{Im}\big(z_1^2\big),
\qquad
\rho_9=\operatorname{Re}\big(z_2^2\big),
\qquad
\rho_{10}=\operatorname{Im}\big(z_2^2\big).
\end{gather*}
The Hilbert series is computed to be{\samepage
\begin{gather*}
\operatorname{Hilb}_{\mathbb{R}[\mathbb{C}^2]^{\mathbb{Z}_2}|\mathbb{R}}(t)=\frac{1+6t^2+t^4}{(1-t^2)^4}=
\sum\limits_{n=0}^\infty{3+2n\choose3}t^{2n}
\\
 \phantom{\operatorname{Hilb}_{\mathbb{R}[\mathbb{C}^2]^{\mathbb{Z}_3}|\mathbb{R}}(t)}
{}=1+10t^2+35t^4+84t^6+165t^8+286t^{10}+\cdots.\tag*{$\triangle$}
\end{gather*}}

$\boldsymbol{N=3.}$ In this case, the Hilbert series is computed to be
\begin{gather*}
\operatorname{Hilb}_{\mathbb{R}[\mathbb{C}^2]^{\mathbb{Z}_3}|\mathbb{R}}(t)
=\frac{1-2t+5t^2-2t^3+t^4}{(1-t)^2(1-t^3)^2}
\\
\phantom{\operatorname{Hilb}_{\mathbb{R}[\mathbb{C}^2]^{\mathbb{Z}_3}|\mathbb{R}}(t)}
{}=1+4t^2+8t^3+9t^4+20t^5+30t^6+36t^7+57t^8+\cdots.
\end{gather*}
It is moreover easy to see that the $4$ quadratic invariants are given by $z_1\cc z_1$, $z_2\cc z_2$, and the real and
imaginary parts of $z_1z_2$. \hfill $\triangle$

\textbf{General $\boldsymbol{N}$.} In general, equation~\eqref{MolSeries} yields the Hilbert series
\begin{gather*}
\operatorname{Hilb}_{\mathbb{R}[\mathbb{C}^2]^{\mathbb{Z}_N}|\mathbb{R}}(t)=
\frac{1}{N}\sum\limits_{g\in\mathbb{Z}_N}\frac{1}{\det(\operatorname{id}-gt)}=\frac{1}{N}\sum\limits_{\zeta^N=
1}\frac{1}{(1-\zeta t)^2(1-\zeta^{-1}t)^2}.
\end{gather*}
Applying Gessel's formula~\cite[Theorem 4.2]{Gessel}, it follows that
$\operatorname{Hilb}_{\mathbb{R}[\mathbb{C}^2]^{\mathbb{Z}_N}|\mathbb{R}}(t)$ is given by the $x^2y^2$-coef\/f\/icient
in the formal power series expansion of
\begin{gather*}
\frac{1}{(1-x)(1-y)-t^2}%\cdot
\left(\frac{x\big(1-x-t^2\big)(1-x)^{N-1}}{(1-x)^N-t^N}+\frac{y(1-y-t^2)(1-y)^{N-1}}{(1-y)^N-t^N}-xy\right).
\end{gather*}
The $x^2y^2$-coef\/f\/icient is
\begin{gather*}
\operatorname{Hilb}_{\mathbb{R}[\mathbb{C}^2]^{\mathbb{Z}_N}|\mathbb{R}}(t)=
\frac{1+t^2+2Nt^N-t^{2N}-2Nt^{N+2}-t^{2N+2}}{(1-t^2)^3(1-t^N)^2}.
\end{gather*}
From this expression, it is easy to compute that if $N>2$, then the $t^2$-coef\/f\/icient of the Hilbert series is
$4$.
Similarly, the $t^3$-coef\/f\/icient vanishes for $N>3$.

\subsubsection{Binary dihedral groups} Suppose now that $G$ is a~binary dihedral group $\mathbb{D}_N$ of order $4N$ for
$N\geq1$, which in complex coordinates is generated by the two elements,
$\operatorname{diag}(\omega_{2N},\omega_{2N}^{-1})$ and
\begin{gather}
\label{eq-Defb}
b=\left[
\begin{matrix}
0&1
\\
-1&0
\end{matrix}
\right],
\end{gather}
where $\omega_{2N}$ is a~primitive $2N$th root of unity.
The action of these generators on $\mathbb C^2\times\overline{\mathbb C^2}$ is given by the matrices
$\alpha_{2N}=\operatorname{diag}\big(\omega_{2N},\omega_{2N}^{-1},\omega_{2N}^{-1},\omega_{2N}\big)$
and $\beta=\operatorname{diag}(b,b)$, the latter in $2\times2$ blocks.

To apply equation~\eqref{MolSeries}, note that the $4N$ elements of $\mathbb{D}_N$ are given by the $2N$ powers of
$\alpha_{2N}$ as well as elements of the form $\alpha^k\beta$ with $0\leq k\leq2N-1$.
A simple computation demonstrates that $\det(\operatorname{id}-\alpha^k\beta t)=(1+t^2)^2$ does not depend on $k$ so
that the sum in equation~\eqref{MolSeries} can be split into a~sum over the cyclic group $\langle\alpha\rangle$ and its
complement in $\mathbb{D}_N$.
That is,
\begin{gather*}
\operatorname{Hilb}_{\mathbb{R}[\mathbb{C}^2]^{\mathbb{D}_N}|\mathbb{R}}(t)=
\frac{1}{4N}\sum\limits_{g\in\mathbb{D}_N}\frac{1}{\det(\operatorname{id}-gt)}
=\frac{1}{4N}\left(\sum\limits_{k=1}^{2N}\frac{1}{\det(\operatorname{id}-\alpha^k t)}+\sum\limits_{k=
1}^{2N}\frac{1}{\det(\operatorname{id}-\alpha^k\beta t)}\right)
\\
\phantom{\operatorname{Hilb}_{\mathbb{R}[\mathbb{C}^2]^{\mathbb{D}_N}|\mathbb{R}}(t)}
{}=\frac{1}{4N}\left(\sum\limits_{\zeta^{2N}=1}\frac{1}{(1-\zeta t)^2(1-\zeta^{-1}t)^2}\right)+\frac{1}{2(1+t^2)^2},
\end{gather*}
where the f\/inal sum corresponds to the case of a~cyclic group treated above.
It follows again by an application of Gessel's formula that the Hilbert series is given by
$\operatorname{Hilb}_{\mathbb{R}[\mathbb{C}^2]^{\mathbb{D}_N}|\mathbb{R}}(t)=Q(t)/P(t)$ where
\begin{gather*}
Q(t)=1+3t^4+(2N-1)t^{2N}+(2N+3)t^{2N+2}-(2N+3)t^{2N+4}
\\
\phantom{Q(t)=}
{}-(2N-1)t^{2N+6}-3t^{4N+2}-t^{4N+6},
\\
P(t)=\big(1-t^2\big)\big(1-t^4\big)^2\big(1-t^{2N}\big)^2.
\end{gather*}
In particular, as $Q(t)/P(t)$ are even functions, the Hilbert series coef\/f\/icients vanish in odd orders.
Moreover, one computes that the $t^2$-coef\/f\/icient is $4$ if $N=1$ and $1$ if $N>1$.

\subsubsection{Binary tetrahedral, octahedral, and icosahedral groups}

The three remaining subgroups of $\operatorname{SU}_2$ are the binary tetrahedral group $\mathbb{T}_{24}$, the binary
octahedral group $\mathbb{O}_{48}$, and the binary icosahedral group $\mathbb{I}_{120}$.
For our purposes, it will be suf\/f\/icient to note that $\mathbb{T}_{24}$ and $\mathbb{O}_{48}$ both contain
a~subgroup isomorphic to $\mathbb{Z}_4$ which in complex coordinates is generated by
$\operatorname{diag}(\sqrt{-1},-\sqrt{-1})$, and hence coincides with $\mathbb{Z}_4$ above.
Similarly, $\mathbb{I}_{120}$ contains a~subgroup isomorphic to $\mathbb{Z}_{10}$ which in complex coordinates is
generated by $\operatorname{diag}(-\omega_5^3,-\omega_5^2)$, and hence is conjugate to the action of $\mathbb{Z}_{10}$
given above.

\subsection[${[}-1,1,1{]}$ and ${[}-1,1,2{]}$ are not orbifolds]{$\boldsymbol{[-1,1,1]}$ and $\boldsymbol{[-1,1,2]}$ are not orbifolds}
\label{finalargument}

Our f\/inal aim is to show that $\operatorname{Hilb}_{\mathbb{R}[M_A]|\mathbb{R}}(t)$ for $A=[-1,1,1]$ or $[-1,1,2]$
cannot coincide with $\operatorname{Hilb}_{\mathbb{R}[\mathbb{C}^2]^G|\mathbb{R}}(t)$ for any f\/inite subgroup
$G<\operatorname{U}_2$.
The argument will be based on the following observation.
Suppose $H\leq G$ and consider
\begin{gather*}
\operatorname{Hilb}_{\mathbb{R}[\mathbb{C}^2]^H|\mathbb{R}}(t)=\sum\limits_{k=0}^\infty a_k t^k,
\qquad
\operatorname{Hilb}_{\mathbb{R}[\mathbb{C}^2]^G|\mathbb{R}}(t)=\sum\limits_{k=0}^\infty b_k t^k.
\end{gather*}
Then any polynomial invariant under $G$ is also invariant under $H$, implying that $b_k\leq a_k$ for each~$k$.
Moreover, note that any f\/inite subgroup of $\operatorname{U}_1$ acting on $\mathbb{C}^2$ as scalar multiplication is
cyclic, and in $(z_1,z_2,\cc z_1,\cc z_2)$-coordinates is generated by
$\operatorname{diag}(\omega_N,\omega_N,\omega_N^{-1},\omega_N^{-1})$ for some $N$.
Permuting coordinates, this action is conjugate to the action generated by
$\operatorname{diag}(\omega_N,\omega_N^{-1},\omega_N,\omega_N^{-1})$ so that the algebra of real invariant polynomials
is isomorphic to the algebra of invariants of a~cyclic subgroup of $\operatorname{SU}_2$, see
Subsection~\ref{subsubsec:CyclicGps}.

We use the notation $(L/L_K;R/R_K)$ to indicate a~f\/inite subgroup of $\operatorname{U}_2$, where $L_K\unlhd L$ are
f\/inite subgroups of $\operatorname{U}_1$, $R_K\unlhd R$ are f\/inite subgroups of $\operatorname{SU}_2$, and
$(L/L_K;R/R_K)$ contains both~$L_K$ and~$R_K$ as subgroups.
Note that $\unlhd$ indicates a~normal subgroup.

\subsubsection[${[}-1,1,1{]}$ is not an orbifold]{$\boldsymbol{[-1,1,1]}$ is not an orbifold}

Let $M_0$ denote the reduced space associated to the weight matrix
$[-1,1,1]$.
Recall (see equation~\eqref{Hilb111}) that in this case, $\operatorname{Hilb}_{\mathbb R[M_0]|\mathbb
R}(t)=1+8t^2+27t^4+\cdots$.
Based on the $t^2$-coef\/f\/icient, it follows that if $\operatorname{Hilb}_{\mathbb R[M_0]|\mathbb
R}(t)=\operatorname{Hilb}_{\mathbb{R}[\mathbb{C}^2]^{G}|\mathbb{R}}(t)$ for some f\/inite $G<\operatorname{U}_2$, then
$G$ contains no cyclic subgroups of order greater than $2$ and no binary dihedral subgroups.
This eliminates all of the nine families of f\/inite subgroups of $\operatorname{U}_2$ other than Type~1 and Type~3, see Appendix~\ref{CoxApp}.

Suppose $G$ is a~Type 1 group of the form $(\mathbb{Z}_{2m}/\mathbb{Z}_f;\mathbb{Z}_{2n}/\mathbb{Z}_g)_d$, and then as
$G$ contains $\mathbb{Z}_f<\operatorname{U}_2$ and $\mathbb{Z}_g<\operatorname{SU}_2$, we have that $f\leq2$ and
$g\leq2$.
As $f\equiv g\mod2$, there are only two cases.

\textbf{Type 1, $\boldsymbol{f=g=2}$.} In this case, $G$ is of the form
$(\mathbb{Z}_{2r}/\mathbb{Z}_2;\mathbb{Z}_{2r}/\mathbb{Z}_2)_d$ where $r$ is a~positive integer and $d\leq r$ is
relatively prime to $r$.
If $r=1$, then $d=1$ and $G=\mathbb{Z}_2<\operatorname{SU}_2$, whose Molien series was computed in
Subsection~\ref{Coxeter} and does not coincide with $\operatorname{Hilb}_{\mathbb R[M_0]|\mathbb R}(t)$.
So assume $r\geq2$.
Then $G$ contains a~subgroup generated by
$\omega_{2r}\operatorname{diag}(\omega_{2r}^d,\omega_{2r}^{-d})
=\operatorname{diag}(\omega_{2r}^{d+1},\omega_{2r}^{1-d})$
for a~primitive $2r$th root of unity $\omega_{2r}$,
which in $(z_1,z_2,\cc z_1,\cc z_2)$-coordinates is given by the $4\times4$ matrix
$\operatorname{diag}(\omega_{2r}^{d+1},\omega_{2r}^{1-d},\omega_{2r}^{-d-1},\omega_{2r}^{d-1})$.

To show that the dimension of quadratic invariants f\/ixed by $\langle\alpha\rangle$ has dimension strictly less than
$8$, and hence that no group containing $\langle\alpha\rangle$ can have the same Hilbert series as
$\operatorname{Hilb}_{\mathbb R[M_0]|\mathbb R}(t)$, consider the action of $\alpha$ on the quadratic polynomials in
$z_1$, $z_2$, $\cc z_2$, $\cc z_2$.
With respect to the basis
\begin{gather*}
\big\{
z_1\cc z_1,\,
z_2\cc z_2,\,
z_1z_2,\,
\cc z_1\cc z_2,\,
z_1\cc z_2,\,
\cc z_1z_2,\,
z_1^2,\,
\cc z_1^2,\,
z_2^2,\,
\cc z_2^2
\big\}
\end{gather*}
for the quadratic polynomials, $\alpha$ has $10\times10$ matrix
\begin{gather*}
Q_\alpha=
\operatorname{diag}
\big(1,\,1,\,\omega_{2r}^2,\,\omega_{2r}^{-2},\,\omega_{2r}^{2d},\,\omega_{2r}^{-2d},\,
\omega_{2r}^{2d+2},\,\omega_{2r}^{-2d-2},\,\omega_{2r}^{-2d+2},\,\omega_{2r}^{2d-2}\big).
\end{gather*}
Using the trace formula~\cite[Lemma 2.2.2]{SturmfelsBook} used to prove Molien's formula, the dimension of the
quadratic polynomials invariant under the action of $\langle\alpha\rangle$ is
\begin{gather*}
\frac{1}{|Q_\alpha|}\sum\limits_{k=1}^{|Q_\alpha|}\operatorname{trace}{Q_\alpha^k},
\end{gather*}
where $|Q_\alpha|$ denotes the the order of $Q_\alpha$.
Note that $|Q_\alpha|$ clearly divides $r$ and moreover that the above formula holds if $|Q_\alpha|$ is replaced by any
positive multiple of $|Q_\alpha|$.
Clearly, $\omega_{2r}^2$~and~$\omega_{2r}^{-2}$ are primitive $r$th roots of unity; similarly, as $d$ is relatively
prime to~$r$, $\omega_{2r}^{2d}$ and $\omega_{2r}^{-2d}$ are primitive $r$th roots of unit as well.
With this, as the sum of all $r$th roots of unity is zero, we have that the dimension of quadratic
$\langle\alpha\rangle$-invariants is
\begin{gather*}
\frac{1}{r}\sum\limits_{k=
1}^{r}2+\omega_{2r}^{(2d+2)k}+\omega_{2r}^{-(2d+2)k}+\omega_{2r}^{(-2d+2)k}+\omega_{2r}^{-(2d+2)k}\leq\frac{1}{r}(6r)=6.
\end{gather*}
Therefore, $G$ cannot contain $\alpha$.
\hfill $\triangle$

 \textbf{Type 1, $\boldsymbol{f=g=1}$.} If $G=(\mathbb{Z}_r/1;\mathbb{Z}_r/1)_d$ for $r$ even and $d<r$
relatively prime to $r$, then $G$ contains
$\alpha=\omega\operatorname{diag}(\omega_r^d,\omega_r^{-d})=\operatorname{diag}(\omega_r^{d+1},\omega_r^{-d+1})$.
If $r=2$, then $G$ is trivial, so assume $r\geq4$.
The action on quadratic polynomials is given by
\begin{gather*}
Q_\alpha=
\operatorname{diag}
\big(1,\,1,\,\omega_r^2,\,\omega_r^{-2},\,\omega_r^{2d},\,\omega_r^{-2d},\,
\omega_r^{2d+2},\,\omega_r^{-2d-2},\,\omega_r^{-2d+2},\,\omega_r^{2d-2}\big).
\end{gather*}
As $\omega_r^2$, $\omega_r^{-2}$, $\omega_r^{2d}$, and $\omega_r^{-2d}$ are primitive $r/2$nd roots of unity, the
dimension of quadratic $\langle\alpha\rangle$-invariants is
\begin{gather*}
\frac{1}{r}\sum\limits_{k=
1}^{r}2+\omega_r^{(2d+2)k}+\omega_r^{-(2d+2)k}+\omega_r^{(-2d+2)k}+\omega_r^{-(2d+2)k}\leq\frac{1}{r}(6r)=6.
\end{gather*}
Again, $G$ cannot contain $\alpha$.
It follows that $G$ cannot be a~Type 1 group.
\hfill $\triangle$

 \textbf{Type 3.} Suppose $G$ is a~Type 3 group of the form either
$(\mathbb{Z}_{4m}/\mathbb{Z}_{2m};\mathbb{D}_l/\mathbb{Z}_{2l})$ or
$(\mathbb{Z}_{4m}/\mathbb{Z}_m;$ $\mathbb{D}_l/\mathbb{Z}_l)$ with $m$ and $l$ odd.
As $G$ cannot contain cyclic subgroups of $\operatorname{U}_2$ or $\operatorname{SU}_2$ of orders larger than $2$, we
need only consider the case of $m=l=1$.
Both $(\mathbb{Z}_4/\mathbb{Z}_2;\mathbb{D}_1/\mathbb{Z}_2)$ and $(\mathbb{Z}_4/1;\mathbb{D}_1/1)$ contain the element
$\sqrt{-1}b$ where $b$ the $2\times2$ matrix def\/ined in equation~\eqref{eq-Defb} above.
One computes with Singular that the Hilbert series of invariants for the group generated by this element is
\begin{gather}
\label{MolZ4D1}
\frac{1+t^2}{(1-t^2)^2(1-t)^2}=1+2t+6t^2+10t^3+19t^4+\cdots.
\end{gather}
Hence, this cannot occur as an element of $G$.
\hfill $\triangle$

With this, it follows that no such $G$ exists, and $[-1,1,1]$ does not admit a~$\mathbb{Z}$-graded regular
symplectomorphism with an orbifold.

\subsubsection[${[}-1,1,2{]}$ is not an orbifold]{$\boldsymbol{[-1,1,2]}$ is not an orbifold}

Let $M_0$ denote the reduced space associated to the weight matrix
$[-1,1,2]$.
Recall (equation~\eqref{Hilb111}) that in this case, $\operatorname{Hilb}_{\mathbb R[M_0]|\mathbb
R}(t)=1+4t^2+6t^3+9t^4+16t^5+\cdots$.
Based on the $t^3$-coef\/f\/icient, it follows that if $\operatorname{Hilb}_{\mathbb R[M_0]|\mathbb
R}(t)=\operatorname{Hilb}_{\mathbb{C}[z_1,\cc z_1,z_2,\cc z_2]^G}(t)$ for some f\/inite $G<\operatorname{U}_2$, then
$G$ contains no cyclic subgroups of order other than $3$ and no binary dihedral subgroups.
This eliminates all of the nine families of f\/inite subgroups of $\operatorname{U}_2$ other than Type 1 and Type 3.

Suppose $G$ is a~Type 1 subgroup of the form $(\mathbb{Z}_{2m}/\mathbb{Z}_f;\mathbb{Z}_{2n}/\mathbb{Z}_g)_d$, and then
as $G$ contains $\mathbb{Z}_f<\operatorname{U}_2$ and $\mathbb{Z}_g<\operatorname{SU}_2$, we have that $f=1$ or $3$ and
$g=1$ or $3$.

\textbf{Type 1, $\boldsymbol{f=g=3}$.} In this case, $G$ is of the form
$(\mathbb{Z}_{3r}/\mathbb{Z}_3;\mathbb{Z}_{3r}/\mathbb{Z}_3)_d$ where $r$ is a~positive even integer and $d<r$ is
relatively prime to $r$.
Then $G$ contains $\mathbb{Z}_3<\operatorname{SU}_2$ as well as the subgroup of $\operatorname{U}_1$ generated by
scalar multiplication by a~primitive $3$rd root of unity $\omega_3$.
As the quadratic invariants of $\mathbb{Z}_3<\operatorname{SU}_2$ are $z_1\cc z_1$, $z_2\cc z_2$, $z_1z_2$, and~$\cc
z_1\cc z_2$, we need only note that~$z_1z_2$ is not invariant under scalar multiplication by~$\omega_3$, so that the
space of quadratic $G$-invariants has dimension strictly less than $4$.
It follows that the Hilbert series of $G$-invariants cannot coincide with $\operatorname{Hilb}_{\mathbb R[M_0]|\mathbb
R}(t)$, and $G$ cannot be of this type.
\hfill $\triangle$

 \textbf{Type 1, $\boldsymbol{f=3}$, $\boldsymbol{g=1}$.} In this case,
$G=(\mathbb{Z}_{3r}/\mathbb{Z}_3;\mathbb{Z}_r/1)_d$ where $r$ is a~positive even integer and $d<r$ is relatively prime
to $r$.
Note that $\mathbb{Z}_3<\operatorname{U}_1$ is a~proper subgroup of $G$.
Then
$G\ni\omega_{3r}\operatorname{diag}(\omega_{3r}^{3d},\omega_{3r}^{-3d})
=\operatorname{diag}(\omega_{3r}^{3d+1},\omega_{3r}^{-3d+1})$
for a~primitive $3r$th root of unity $\omega_{3r}$, which is given in $(z_1,z_2,\cc z_1,\cc z_2)$-coordinates
by $\operatorname{diag}(\omega_{3r}^{3d+1},\omega_{3r}^{-3d+1},\omega_{3r}^{-3d-1},\omega_{3r}^{3d-1})$.
The quadratic invariants of the action of $\mathbb{Z}_3<\operatorname{U}_1$ are $z_1\cc z_1$, $z_2\cc z_2$, $z_1\cc
z_2$, and $\cc z_1z_2$.
It is easy to see that $z_1\cc z_2$ is not invariant under the action of $\alpha$ for $r>2$, implying that the space of
quadratic $G$-invariants has dimension strictly less than $4$.
So assume $r=2$, and then the Hilbert series of $G=(\mathbb{Z}_6/\mathbb{Z}_3;\mathbb{Z}_2/1)_1$ is computed on
Singular to be
\begin{gather*}
\frac{1-2t+5t^2-2t^3+t^4}{(1-t^3)^2(1-t)^2}=1+2t^2+4t^3+3t^4+8t^5+{12}+\cdots,
\end{gather*}
which does not coincide $\operatorname{Hilb}_{\mathbb R[M_0]|\mathbb R}(t)$.
Hence $G$ cannot be of this type.
\hfill $\triangle$

 \textbf{Type 1, $\boldsymbol{f=1}$, $\boldsymbol{g=3}$.} In this case, $G$ is of the form
$(\mathbb{Z}_r/1;\mathbb{Z}_{3r}/\mathbb{Z}_3)_d$ where $r$ is a~positive even integer and $d<r$ is relatively prime to
$r$.
Then $G$ contains a~subgroup generated by
$\omega_{3r}^3\operatorname{diag}(\omega_{3r}^d,\omega_{3r}^{-d})
=\operatorname{diag}(\omega_{3r}^{d+3},\omega_{3r}^{-d+3})$
for a~primitive $3r$th root of unity $\omega_{3r}$.
Note that $G$ contains $\mathbb{Z}_3<\operatorname{SU}_2$, whose quadratic invariants are again spanned by $z_1\cc
z_1$, $z_2\cc z_2$, $z_1z_2$, and $\cc z_1\cc z_2$.
If $r>2$, then $z_1z_2$ is not invariant under the above, implying that the quadratic invariants have dimension
strictly less than $4$.
If $r=2$, then $G=(\mathbb{Z}_2/1;\mathbb{Z}_6/\mathbb{Z}_3)_1$ is in fact a~subgroup of $\operatorname{SU}_2$
isomorphic to $\mathbb{Z}_3$.
We again conclude that $G$ cannot be of this type.
\hfill $\triangle$

\textbf{Type 1, $\boldsymbol{f=g=1}$.} If $G=(\mathbb{Z}_r/1;\mathbb{Z}_r/1)_d$ for even $r$ and $d<r$
relatively prime to~$r$, then~$G$ is generated by
$\alpha=\omega_r\operatorname{diag}(\omega_r^d,\omega_r^{-d})=\operatorname{diag}(\omega_r^{d+1},\omega_r^{-d+1})$ If
$r=2$, then~$G$ is trivial, so assume~$r\geq4$.
If $d=1$, then $\alpha=\operatorname{diag}(\omega_r^2,1)$, and if $d=r-1$, then
$\alpha=\operatorname{diag}(1,\omega_r^2)$.
In either case, $G$~has nontrivial linear invariants, so that as $\operatorname{Hilb}_{\mathbb R[M_0]|\mathbb R}(t)$
has zero $t$-coef\/f\/icient, we may exclude these cases.
So assume $1<d<r-1$, and then as~$d$ must be relatively prime to~$r$, it must be that $r\geq8$.

Now, with respect to the following basis for the cubic monomials,
\begin{gather*}
\big\{
z_1^3,\,\cc z_1^3,\,z_2^3,\,\cc z_2^3,\,z_1^2z_2,\,\cc z_1^2\cc z_2,\,z_1z_2^2,\,
\cc z_1\cc z_2^2,\,z_1^2\cc z_2,\,\cc z_1^2z_2,\,z_1\cc z_2^2,
\\
\phantom{\big\{}
\cc z_1z_2^2,\,z_1^2\cc z_1,\,z_1\cc z_1^2,\,z_2^2\cc z_2,\,z_2\cc z_2^2,\,z_1\cc z_1z_2,\,
z_1\cc z_1\cc z_2,\,z_1z_2\cc z_2,\,\cc z_1z_2\cc z_2
\big\}
\end{gather*}
the action of $\alpha$ is given by
\begin{gather*}
Q_\alpha=
\operatorname{diag}\big(\omega_r^{3d+3},\,\omega_r^{-3d-3},\,\omega_r^{-3d+3},\,
\omega_r^{3d-3},\,\omega_r^{d+3},\,\omega_r^{-d-3},\,\omega_r^{-d+3},\,\omega_r^{d-3},\,\omega_r^{3d+1},\,\omega_r^{-3d-1},
\\
\phantom{Q_\alpha=\operatorname{diag}\big(}
\omega_r^{3d-1},\,\omega_r^{-3d+1},\,\omega_r^{d+1},\,\omega_r^{-d-1},\,
\omega_r^{-d+1},\,\omega_r^{d-1},\,\omega_r^{-d+1},\,\omega_r^{d-1},\,\omega_r^{d+1},\,\omega_r^{-d-1}\big).
\end{gather*}
Applying the trace formula~\cite[Lemma~2.2.2]{SturmfelsBook}, we have that the dimension of the cubic polynomials
invariant under the action of $\langle\alpha\rangle$ is given by
\begin{gather}
\nonumber
\frac{1}{r}\sum\limits_{k=
1}^{r}\omega_r^{(3d+3)k}+\omega_r^{(-3d-3)k}+\omega_r^{(-3d+3)k}+\omega_r^{(3d-3)k}+\omega_r^{(d+3)k}+\omega_r^{(-d-3)k}
+\omega_r^{(-d+3)k}\\
\qquad
{}+\omega_r^{(d-3)k}+\omega_r^{(3d+1)k}+\omega_r^{(-3d-1)k}+\omega_r^{(3d-1)k}
+\omega_r^{(-3d+1)k}+\omega_r^{(d+1)k}+\omega_r^{(-d-1)k}\nonumber\\
\qquad {}+\omega_r^{(-d+1)k}+\omega_r^{(d-1)k}+\omega_r^{(-d+1)k}
+\omega_r^{(d-1)k}+\omega_r^{(d+1)k}+\omega_r^{(-d-1)k}.
\label{CubicTr}
\end{gather}
For $k=1$, each term in the above sum is a~primitive $s$th root of unity for some $s$ that divides~$r$ so that unless
a~term is equal to~$1$, the sum over~$k$ of that term vanishes.
Recalling that $1<d<r-1$ and $r\geq8$, it is clear that $\omega_r^{\pm d\pm1}\neq1$, so that the sum of each
$\omega_r^{(\pm d\pm1)k}$ vanishes.
Similarly, as~$d$ is relatively prime to~$r$ and hence invertible mod~$r$, it is easy to see that at most two of the
following congruences can be true mod~$r$:
\begin{alignat*}{4}
& 3d+3\equiv0,\qquad && 3d-3\equiv0,\qquad && d+3\equiv0, &
\\
& d-3\equiv0,\qquad && 3d+1\equiv0,\qquad && 3d-1\equiv0. &
\end{alignat*}
Therefore, when $k=1$, at most four of the terms in equation~\eqref{CubicTr} can be equal to~$1$, and the dimension of
cubic invariants is bounded by $\frac{1}{r}(4r)=4$.
As the dimension of cubic invariants on $M_0$ is six, we have excluded all groups in this case.
\hfill $\triangle$

\textbf{Type 3, $\boldsymbol{(\mathbb{Z}_{12}/\mathbb{Z}_3;\mathbb{D}_3/\mathbb{Z}_3)}$.} As in the
case of a~Type 1 group with $f=g=3$, this group contains $\mathbb{Z}_3<\operatorname{SU}_2$ as well as the subgroup of
$\operatorname{U}_1$ generated by scalar multiplication by a~primitive $3$rd root of unity $\omega_3$; hence the
Hilbert series of $G$-invariants cannot coincide with $\operatorname{Hilb}_{\mathbb R[M_0]|\mathbb R}(t)$.~\hfill$\triangle$

\textbf{Type 3, $\boldsymbol{(\mathbb{Z}_4/1;\mathbb{D}_3/\mathbb{Z}_3)}$.} This group has six
elements and is generated by $\alpha=\operatorname{diag}(\omega_3,\omega_3^2)$ and $\beta=\sqrt{-1}b$, where $b$ is
def\/ined in equation~\eqref{eq-Defb}.
The Hilbert series is given by
\begin{gather*}
\operatorname{Hilb}_{\mathbb{R}[\mathbb{C}^2]^G|\mathbb{R}}(t)=\frac{1+t^2+2t^3+t^4+t^6}{(1-t^2)^2(1-t^3)^2}=
1+3t^2+4t^3+6t^4+\cdots,
\end{gather*}
which does not coincide with $\operatorname{Hilb}_{\mathbb R[M_0]|\mathbb R}(t)$.
\hfill $\triangle$

\textbf{Type 3, $\boldsymbol{(\mathbb{Z}_{12}/\mathbb{Z}_3;\mathbb{D}_1/1)}$.} This group has six
elements and is generated by $\omega_{12}b$.
The Hilbert series is given by
\begin{gather*}
\operatorname{Hilb}_{\mathbb{R}[\mathbb{C}^2]^G|\mathbb{R}}(t)=
\frac{1-t+t^2+2t^3+2t^5+t^6-t^7+t^8}{(1-t^6)(1-t^3)(1-t^2)(1-t)}=1+2t^2+4t^3+\cdots,
\end{gather*}
which does not coincide with $\operatorname{Hilb}_{\mathbb R[M_0]|\mathbb R}(t)$.
\hfill $\triangle$

\textbf{Type 3, $\boldsymbol{(\mathbb{Z}_4/1;\mathbb{D}_1/1)}$.} The only nontrivial element of this
group is $\sqrt{-1}b$.
The Hilbert series was computed in equation~\eqref{MolZ4D1} above and does not coincide with
$\operatorname{Hilb}_{\mathbb R[M_0]|\mathbb R}(t)$.
\hfill $\triangle$

\appendix

%\pdfbookmark[1]{Appendix: Finite subgroups of ${\rm U}_2$}{appendix}
\section[Finite subgroups of ${\rm U}_2$]{Finite subgroups of $\boldsymbol{{\rm U}_2}$}
\label{CoxApp}

For the convenience of the reader, we recall the classif\/ication of f\/inite subgroups of $\operatorname{U}_2$ given
by~\cite{Coxeter, DuVal}.
We follow~\cite{FalbelPaupert}; see also~\cite{dunbar}.

For $l\in\operatorname{U}_1$ and $r\in\operatorname{SU}_2$, we let $(l,r)$ denote the element of $\operatorname{U}_2$
given by the scalar multiple~$lr$ of~$r$.
Note that every element of $\operatorname{U}_2$ arises in this way and that the expression is unique up to
$(l,r)=(-l,-r)$.
Let $L_K\unlhd L<\operatorname{U}_1$ and $R_K\unlhd R<\operatorname{SU}_2$ be f\/inite subgroups of~$\operatorname{U}_1$ and~$\operatorname{SU}_2$, respectively, such that~$L/L_K$ is isomorphic to~$R/R_K$, and let
$\phi\colon  L/L_K\to R/R_K$ be an isomorphism.
Then the group $(L/L_K;R/R_K)_\phi$ is def\/ined as{\samepage
\begin{gather*}
(L/L_K;R/R_K)_\phi=\big\{(l,r)\in L\times R:\phi(lL_K)=rR_K\big\}.
\end{gather*}
Note that $\phi$ is omitted if it is obvious, and $(L/L_K;R/R_K)_\phi$ has order $|R||L_K|/2$.}

Let $\mathbb{Z}_k$ denote a~cyclic subgroup of order $k$.
Note that $\mathbb{Z}_k<\operatorname{U}_1$ is generated by a~primitive $k$th root of unity $\omega_k$, while
$\mathbb{Z}_k<\operatorname{SU}_2$ is generated by $\operatorname{diag}(\omega_k,\omega_k^{-1})$.
The distinction will be clear from the context.
Let $\mathbb{D}_p$ denote the binary dihedral group of order $4p$, and let $\mathbb{T}_{24}$, $\mathbb{O}_{48}$, and
$\mathbb{I}_{120}$ denote the binary tetrahedral, octahedral, and icosahedral groups, respectively.

The f\/inite subgroups of $\operatorname{U}_2$ are given by the following list.
\begin{enumerate}[Type 1.]\itemsep=0pt
\item $(\mathbb{Z}_{2m}/\mathbb{Z}_f;\mathbb{Z}_{2n}/\mathbb{Z}_g)_d$ where $f\equiv g\mod2$, and $d$ is
relatively prime to $2m/f=2n/g$ and indicates the isomorphism
$\mathbb{Z}_{2m}/\mathbb{Z}_f\to\mathbb{Z}_{2n}/\mathbb{Z}_g$ sending the class of 1 to the class of $d$, \item
$(\mathbb{Z}_{2m}/\mathbb{Z}_{2m};\mathbb{D}_l/\mathbb{D}_l)$, \item
$(\mathbb{Z}_{4m}/\mathbb{Z}_{2m};\mathbb{D}_l/\mathbb{Z}_{2l})$ and
$(\mathbb{Z}_{4m}/\mathbb{Z}_m;\mathbb{D}_l/\mathbb{Z}_l)$ for $m$ and $l$ odd, \item
$(\mathbb{Z}_{4m}/\mathbb{Z}_{2m};\mathbb{D}_{2l}/\mathbb{D}_l)$, \item
$(\mathbb{Z}_{2m}/\mathbb{Z}_{2m};\mathbb{T}_{24}/\mathbb{T}_{24})$, \item
$(\mathbb{Z}_{6m}/\mathbb{Z}_{2m};\mathbb{T}_{24}/\mathbb{D}_2)$, \item
$(\mathbb{Z}_{2m}/\mathbb{Z}_{2m};\mathbb{O}_{48}/\mathbb{O}_{48})$, \item
$(\mathbb{Z}_{4m}/\mathbb{Z}_{2m};\mathbb{O}_{48}/\mathbb{T}_{24})$, and \item
$(\mathbb{Z}_{2m}/\mathbb{Z}_{2m};\mathbb{I}_{120}/\mathbb{I}_{120})$.
\end{enumerate}

\subsection*{Acknowledgements}

The authors would like to thank Srikanth Iyengar, Luchezar Avramov, Markus Pf\/laum, Johan Martens, Karl-Heinz
Fieseler, Jedrzej Sniatycki, Gerry Schwarz, Johannes Huebschmann, Michael J. Field and Graeme Wilkin for promptly answering questions, stimulating discussions, and moral support.
We would also like to thank the referees for helpful suggestions and comments.

C.F.~would like to thank the University of Florence for hospitality during the completion of this manuscript.
The research of H.-C.~H.~has been supported by the Center for the Quantum Geometry of Moduli spaces which is funded by the Danish National
Research Foundation, and by the Department of Mathematics of the University of Nebraska at Lincoln.
C.S.~received support from the Center for the Quantum Geometry of Moduli spaces, a~Rhodes College Faculty Development
Endowment Grant, and a~grant to Rhodes College from the Andrew W.~Mellon Foundation.

\pdfbookmark[1]{References}{ref}
\LastPageEnding

\end{document}